\documentclass[12pt]{article}
\usepackage{amssymb, amsmath}
\usepackage{amsfonts}
\usepackage{amsthm}
\usepackage{enumerate}

\theoremstyle{definition}
\newtheorem{theorem}{Theorem}[section]

\newtheorem{definition}[theorem]{Definition}

\newcommand{\bysame}{\mbox{\rule{3em}{.4pt}}\,}

\def\dim{\textrm{ dim }}
\def\Ker{\textrm { ker }}
\def\lra{\longrightarrow}
\def\rank{\textrm{ rank }}
\def\R{\mathbb R}

\def\pr{\pi_R}
\def\pl{\pi_L}
\def\hat{\widehat}
\def\tilde{\widetilde}
\def\bfo{\begin{eqnarray*} }
\def\efo{\end{eqnarray*} }
\def\ba{\begin{eqnarray*} }
\def\ea{\end{eqnarray*} }
\def\beq{\begin{eqnarray}}
\def\eeq{\end{eqnarray}}
\def\supp{\hbox{supp}\,}

\def\det {\hbox{det}}
\def\So{\Sigma_1}
\def\Soo{\Sigma_{1,1}}
\def\tC{\tilde{C}}
\def\tS{\tilde\Sigma}
\def\tht{\theta_3}
\def\L{\Lambda}
\def\ni{\noindent}
\def\tf{\tilde{f}}

\begin{document}

\title{Fourier integral operators with open umbrellas and  seismic inversion for cusp caustics}

\author{Raluca Felea and Allan Greenleaf}
%%\footnote{The second author was partially supported by NSF grants
%%DMS-0551894 and -0853892.}

%%\thanks{The second author was partially supported by NSF grants
%%DMS-0551894 and -0853892.}

%\address{School of Mathematical Sciences \\ Rochester  Institute of Technology\\
%Rochester, NY 14623}
 % \email{rxfsma@rit.edu}
%  %%optional:  \curaddr{current address}%%
%   %%optional:  \urladdr{website address}%%

  %  \address{Department of Mathematics \\ University of Rochester\\ Rochester, NY 14627}
 % \email{allan@math.rochester.edu}
%  %%optional:  \curaddr{current address}%%
%   %%optional:  \urladdr{website address}%%

\date{}

\maketitle

\begin{abstract} In general the composition of Fourier integral operators (FIOs) need not be an FIO. Motivated by  the problem of linearized seismic inversion in the presence of cusp caustics for the background sound speed, we consider  FIOs whose canonical relations have certain two-sided cusp degeneracies, and show that the resulting compositions have wave-front relations in the union of the diagonal and an open umbrella, the simplest type of singular Lagrangian manifold.

\end{abstract}

%\vspace*{1 cm}

\section{Introduction}\label{sec intro}

A fundamental problem concerning Fourier integral operators (FIOs) is that, outside of the standard transverse \cite{
H-book} and clean intersection \cite{DGu,W} calculus, a composition of two FIOs is typically not another FIO. Describing the
operators resulting from the composition of completely general FIOs and placing them in a usable class, with a symbol
calculus, Sobolev space estimates, and the possibility of constructing parametrices under suitable ellipticity conditions, is
at this point a distant prospect. Some progress has been made for specific geometries arising in integral geometry
\cite{GU1} and inverse problems \cite{F,FG,Mar,N}. In all of these works, the compositions were shown to belong to an
existing class, the \emph{pseudodifferential operators with singular symbols}, which are not FIOS and were introduced  to
construct parametrices for operators of real or complex principal type \cite{MU,GuU,Men}. Such operators have  wave front relation, i.e., the wave front
set of the Schwartz kernel, contained in the union of the diagonal and another smooth canonical relation which
intersects the diagonal cleanly.

In the present paper, motivated by a linearized inverse problem from seismology for the acoustic wave equation, we analyze  FIOs having a
special type of  degeneracy, which we call a \emph{flat two-sided cusp}, and show that their composition results in operators which are fundamentally more singular than those in works referred to above.

FIOs with  \emph{two-sided cusps}, i.e., those for which the projections from the
canonical relation both to the left and the right are  cusps, arise naturally when considering generalized Radon transforms. They have previously been studied  in terms of their $L^2$ Sobolev regularity
properties; it was shown in  \cite{CoCu,GS1} that there is a loss of $1/4$ derivative compared with the nondegenerate case of
local canonical graphs.  $L^2\lra L^q$ estimates for FIOs with one-sided cusps have also been obtained, with $T^*T$ making an
appearance through Strichartz estimate  type arguments \cite{GS}; however, such compositions  have not been studied in their
own right. We show that  for the subclass of \emph{flat} two-sided cusps which we formulate,  the normal operator $F^*F$
is not a standard pseudodifferential operator with singular symbol, as it is in the case of fold caustics \cite{N,F,FG}.  This includes the canonical relations underlying the linearized forward
scattering operators $F$   in the seismology problem in the presence of caustics of \emph{cusp} type. While
there is still a large part of its wave front relation, $WF(K_{F^*F})'$, which is contained in the diagonal, $\Delta$, the
remaining portion, $\tC$, is now no longer a \emph{smooth} canonical relation. Rather, $\tC$ has the structure of the
simplest kind of singular Lagrangian manifold, an \emph{ open umbrella} \cite{Gi}.  The diagonal and this open umbrella
intersect in codimension one, and we quantize the variety $\Delta\cup\tC$ by associating to it a class of generalized Fourier
integral operators, $I^m(\Delta\cup\tC)$. This is   in the spirit of the paired Lagrangian  (or $I^{p,l}$) distributions associated to a pair of cleaning
intersecting smooth Lagrangians;  however,  rather than
 nondegenerate phase functions and product type symbols as in  \cite{MU,GuU,Men}, we  use a combination of a degenerate phase function  (whose gradient exhibits normal crossings) and standard symbols. We remark that operators with wave front relation in the union of $\Delta$  and a canonical relation having a conical singularity were considered by Melrose and Uhlmann \cite{MU2}.

As applied to the seismic imaging problem, the results here have negative implications (as do the results of \cite{N,F} in
the case of fold caustics): even microlocally away from $\Delta\cap\tC$, the non-pseudodifferential part of the normal operator has the same order as the
pseudodifferential part, resulting in a strong, nonremovable artifact. However, our work  provides a precise description of the imaging artifacts
resulting from cusp caustics, observed by Nolan and Symes \cite{NS-fold} in a model 2D case.  It would be of interest to obtain Sobolev estimates for  operators associated with $\Delta\cup\tilde{C}$, similar to the results  in \cite{FGP} for certain $I^{p,l}$ classes, including those arising from fold caustics.

We begin in Sec. \ref{sec examps} with a simple example from harmonic analysis   of a FIO having two-sided cusps of the type
that is relevant for the seismology problem, recall some basic singularity theory, and examine the corresponding normal
operator. In Sec. \ref{sec  caustics}, we describe  the linearized  inverse problem from seismology, and show that the
canonical relation underlying the linearized forward scattering operator $F$ is  (i) associated with a two-sided
cusp, but  with the additional (highly nongeneric) features that \linebreak(ii) the cusp points for the left and right projections are
equal, and (iii) the images of the cusp points on both the left and right are involutive (coisotropic). In Sec. \ref{sec
cusps} we define a general class of canonical relations having this structure,
 and then derive a weak normal form, very close to the model examined in Sec.
\ref{sec examps}, for a general flat two-sided cusp. Finally, in Sec. \ref{sec comp}, we use the oscillatory representations
found in Sec. \ref{sec cusps} to analyze the composition $B^*A$ for two FIOs associated with such a canonical relation. We point out that the
technique of deriving weak normal forms for degenerate FIOs and using  them to prove composition theorems has its origins in
\cite{GU1} and was continued in \cite{FG};  different notions of weak normal forms were used to obtain estimates   in
\nolinebreak\cite{GS-fold,GS}.

\section{Example of an FIO with a flat two-sided cusp}\label{sec examps}

\subsection{A generalized Radon transform}\label{subsec gRt}

We start with a simple example of a generalized Radon transform that incorporates  the features of the linearized seismic
inversion problem in which  we are ultimately  interested. In fact, we will later show that perturbations of this example
serve as weak normal forms for the general class of operators with flat two-sided cusps.

The    operator averaging over translates in $\R^4$ of
a curve $\gamma(t)$  such that
\beq\label{eqn rfour}
\dot\gamma,\, \ddot\gamma,\, \dddot\gamma,\, \gamma^{(4)}\hbox{ are linearly independent }
\eeq
is an FIO having a canonical relation which is a
two-sided cusp \cite{GS}, but the same holds for  the curve $\gamma(t)=(t,t^2,t^4)$ in $\R^3$ \cite{GS1,GS-Esc}, and this is the model we will use.

Thus, consider the generalized Radon transform $R_0:\mathcal D'(\R^3)\lra \mathcal D'(\R^3)$,
 \beq\label{eqn Rnot}
 & &R_0f(x)=\int f\left(x-(t,t^2,t^4)\right) \chi(t) dt \\
 &=&\int_{\R^2} e^{i\left(x_2-y_2-(x_1-y_1)^2))\theta_2 + (x_3-y_3-(x_1-y_1)^4)\tht\right)} f(y) \chi(x_1-y_1)\, 1(\theta)
 d\theta_2 d\tht dy,
 \nonumber
 \eeq
 where $\chi\in C_0^\infty(\R)$ is a fixed cutoff function.
The associated canonical relation, the conormal bundle $C_0=N^*Z'\subset \big( T^*\R^3\setminus 0\big)\times \big( T^*\R^3\setminus 0\big)$, where $Z=\supp (K_{R_0})\subset\R^3\times\R^3$ and
$(x,\xi;y,\eta)'=(x,\xi;y,-\eta)$ is   standard notation for the twist map, is
\beq
\label{eqn czero}
C_0&=&\Big\{ \Big(x_1, x_2, x_3, -2(x_1-y_1)\theta_2-4(x_1-y_1)^3\tht,\theta_2, \tht; \nonumber\\
& &\qquad y_1, y_2, y_3, -2(x_1-y_1)\theta_2-4(x_1-y_1)^3\tht,\theta_2, \tht \Big )\, : \\
& &\qquad\quad  x_2-y_2-(x_1-y_1)^2=x_3-y_3-(x_1-y_1)^4=0\Big\}.\nonumber
\eeq
We will show that  $C_0$ has several properties, which remarkably also hold in the totally unrelated seismic imaging problem:
(i) the  projections  both to the left and right, $\pl,\pr:C_0\lra T^*\R^3\setminus 0$, have (Whitney) cusp degeneracies;
(ii) the cusp points for the two projections  are the same, $\Sigma_{1,1}(\pl)=\Sigma_{1,1}(\pr):=\Soo$; and (iii) the
images of the cusp points, $\pl(\Soo)$ and $\pr(\Soo)$, are coisotropic (involutive) submanifolds of $T^*\R^3$. Conditions
(ii) and (iii) are unstable and quite special  among two-sided cusps, i.e., canonical relations satisfying (i); see the Remarks in Sec.
\ref{subsec class} below.

% \bigskip
\subsection{Singularity classes}\label{subsec sing}

We first recall some basic facts about cusps and refer to \cite{Wh,Mo,GoGu} for more details.
 Let  $f : {\R}^N \to {\R}^N$ be a smooth function. We say that $f$ \emph{drops rank simply} at $p$ if rank $(df)_p= N-1$
and if $\left(d(\det\, df)\right)_p\neq 0$, so that   $\So (f) :=\{ x\in {\R}^N: \det (df (x))=0\}$, the corank one critical set
of $f$,  is a smooth hypersurface near $p$.
If $\ker df_p \not\subset T_p\So(f)$, then $f$ has a \emph{fold} singularity, and $f|_{\So(f)}$ is an immersion. Considering
the more degenerate case when $\ker df_p \subset T_p\Sigma_1(f)$, one may choose a nonzero vector field $v$ along $\So(f)$ such
that $v \in\ker df_p$, so that $v$ is tangent to $\So(f)$ at $p$. Let $g$ be a smooth function such that $g|_{\So(f)}=0$ and
$dg_p \neq 0$, e.g., $g=\det\,  df$. Thus,  $dg (v)$ has a zero at $p$.
\begin{definition}  $f$  has a (Whitney) {\em cusp}
at p if $dg (v)$  has a simple zero at $p$.
\end{definition}
$\Soo(f)$, the \emph{cusp set} of $f$, is then codimension 2, and $f|_{\Soo(f)}$ is an immersion.
One may use adapted coordinates to clarify  this. These are local  coordinates  such that $f(x_1, x_2, \dots, x_N)=(x_1,
x_2, \dots, x_{N-1},h(x))$ and $h(0)=0$  \cite{Mo}.  Then, $\So(f)=\{x: \frac{\partial h}{\partial x_N} (0)=0 \}$, and $f$
has a cusp singularity  at $0$ iff $\frac{\partial^2 h}{\partial x_N^2} (0)=0, \ \ \frac{\partial^3 h}{\partial x_N^3} (0)
\neq 0$  and  $\rank [d_x(\frac{\partial h}{\partial x_N}), d_x(\frac{\partial^2 h}{\partial x_N^2}) ]=2$. These conditions
are adapted coordinate-independent,
and the notion of a cusp makes sense for any smooth mapping between $N$-dimensional manifolds.
In suitable smooth coordinates on the domain and range spaces, any map with a cusp singularity can be put into the
local normal form,  $f(x_1,x_2,\dots x_N)=(x_1,x_2, \dots, x_{N-1}, x_{N-1}x_N+\nolinebreak x_N^3)$.

\bigskip

For the canonical relation $C_0$, the projection to the left, $\pi_L:C_0\lra T^*\R^3$, is
$ \pi_L( x_1, x_2, x_3, y_1, \theta_2, \tht) =(x_1, x_2, x_3, -2(x_1-y_1)\theta_2-4(x_1-y_1)^3\tht,\theta_2, \tht)$;
hence, letting $\alpha=2\theta_2+12(x_1-y_1)^2\theta_3$, $\beta= -2(x_1-y_1)$ and $\gamma=-4(x_1-y_1)^3$, one has
$$d\pi_L= \left[ \begin{array} {cccccc} 1 & 0 & 0 & 0 & 0 & 0  \\
0 &  1  & 0 & 0 & 0 & 0  \\
0 & 0 & 1 & 0 & 0 & 0  \\
-\alpha & 0 & 0 &  \alpha & \beta & \gamma \\
0 & 0 & 0 & 0 & 1 & 0\\
0 & 0 & 0 & 0 & 0 & 1\\
\end{array} \right]$$
which has  $\det\, d\pi_L=2\theta_2+12(x_1-y_1)^2\tht$ and  $\ker \pi_L = \R\cdot\frac{\partial}{\partial y_1} $ at  $\Sigma_1(\pl)=\{
2\theta_2+12(x_1-y_1)^2\tht=0 \}$. Since $ \frac{\partial}{\partial y_1}(2\theta_2+12(x_1-y_1)^2\tht )=-24(x_1-y_1)\tht$ and
$\frac{\partial^2}{\partial y_1^2}(2\theta_2+12(x_1-y_1)^2\tht )=24 \tht \neq 0$,  the cusp set is
$$\Soo(\pl)= \{  2\theta_2+12(x_1-y_1)^2\tht= x_1-y_1=0  \} =\{ \theta_2=x_1-y_1=0  \}.$$
Noting that $2\theta_2 + 12 (x_1-y_1)^2\tht$ and $x_1-y_1$ have linearly independent gradients, we see that $\pi_L$ has a
cusp singularity.

Similarly, the projection to the right is
$$ \pi_R(x, y_1, \theta_2, \tht)\!\! =\!\!(y_1, x_2-(x_1-y_1)^2, x_3-(x_1-y_1)^4, -2(x_1-y_1)\theta_2-4(x_1-y_1)^3\tht,\theta_2,
\tht),$$
so that $\det\, d\pi_R=-2\theta_2-12(x_1-y_1)^2\tht$, and  $\ker \pi_R = \frac{\partial}{\partial x_1} $ at
$\So(\pr)=\{2\theta_2+12(x_1-y_1)^2\tht=0\}$. One also has $ \frac{\partial}{\partial x_1}(2\theta_2+12(x_1-y_1)^2\tht
)=24(x_1-y_1)\tht$ and $\frac{\partial^2}{\partial x_1^2}(2\theta_2+12(x_1-y_1)^2\tht )=24 \tht \neq 0$, so that
$\Soo(\pr)=\{ \theta_2=x_1-y_1=0  \}$, and thus $\pi_R$  also has a cusp singularity.

Note that $\Sigma_1(\pr)=\Sigma_1(\pl)$ and $\Sigma_{1,1}(\pr)=\Sigma_{1,1}(\pl)=:\Soo$. The first of these is true for any
canonical relation \cite{H-book}, but the second is a very strong condition; indeed, for general  two-sided cusps, there is no
relationship between $\Soo(\pl)$ and $\Soo(\pr)$.
Furthermore, for the images of these sets, one has
$$ \pi_L (\Sigma_{1})=  \{ \xi_2^3= -\frac{27}{8} \xi_1^2
\xi_3   \},\quad \pi_R (\Sigma_{1})= \{ \eta_2^3= -\frac{27}{8}
  \eta_1^2 \eta_3   \},$$
with $ \pi_L (\Sigma_{1,1})= \{ \xi_1=\xi_2=0 \},\quad \pi_R (\Sigma_{1,1})= \{ \eta_1= \eta_2=0 \}$, resp., their cuspidal
edges.
$ \pi_L (\Sigma_{1,1})$ and $ \pi_R (\Sigma_{1,1})$  are codimension two coisotropic (or involutive)  submanifolds
of $T^*\R^3\setminus 0$ \cite{H-book}; again, this is a nongeneric situation.
  \bigskip

\par Consider next the composition of canonical relations,
$$C_0^t \circ C_0= \big\{(x, \xi; y, \eta):\,  \exists (z,\zeta)\hbox{ s.t. } (x, \xi; z, \zeta) \in C_0^t \hbox{ and } \  (z,
\zeta; y, \eta) \in C_0   \big \}.$$
For such a triple $(x,\xi;z,\zeta;y,\eta)$, one obtains  the following equations:
$$z_2-x_2=(z_1-x_1)^2,\quad z_3-x_3=(z_1-x_1)^4,   $$
$$ \xi_1=\zeta_1=-2(z_1-x_1)\zeta_2-4(z_1-x_1)^3\zeta_3,\, \xi_2=\zeta_2, \, \xi_3=\zeta_3,   $$
$$z_2-y_2=(z_1-y_1)^2,\quad
z_3-y_3=(z_1-y_1)^4,   $$
$$ \eta_1=\zeta_1=-2(z_1-y_1)\zeta_2-4(z_1-y_1)^3\zeta_3,\, \eta_2=\zeta_2, \ \ \eta_3=\zeta_3.$$
Using $2(z_1-x_1)\xi_2+4(z_1-x_1)^3\xi_3 =2(z_1-y_1)\xi_2+4(z_1-y_1)^3\xi_3$,   after simplification  one obtains:
$$(y_1-x_1)\xi_2+2[(z_1-x_1)^3-(z_1-y_1)^3]\xi_3=0, $$
$$(y_1-x_1)[\xi_2 + 2(3z_1^2-3z_1(x_1+y_1)+ x_1^2+y_1^2+x_1y_1)\xi_3]=0  $$
It follows  that the contribution  to  $C_0^t \circ C_0$ from  $\{y_1-x_1=0\}$ is contained in $\Delta$, and  that from
$\{y_1-x_1 \neq 0\}$ is  contained in $\tC_0$, where
\beq\label{eqn tC_0}\nonumber
\tC_0&=&\Big\{ \big(x_1, x_2, x_3, -2(z_1-x_1)\theta_2-4(z_1-x_1)^3\theta_3, \theta_2, \theta_3;\nonumber\\
& & \quad y_1, y_2, y_3, -2(z_1-x_1)\theta_2-4(z_1-x_1)^3\theta_3,\  \theta_2, \ \theta_3\big) : \\
& & \qquad\qquad x\in\R^3,\, y_1,z_1\in\R,\, \theta_3\in\R\backslash 0,\, y_2=  x_2+(y_1-x_1)(2z_1-x_1-y_1),\nonumber\ \\
& & \qquad\qquad y_3=x_3+(y_1-x_1)(2z_1-x_1-y_1)((z_1-x_1)^2+(z_1-y_1)^2), \nonumber\\
& & \qquad\qquad  \theta_2=-2(3z_1^2-3z_1(x_1+y_1)+ x_1^2+y_1^2+x_1y_1)\theta_3 \quad \Big\}\nonumber
\eeq
Notice that $\xi_1=\eta_1= 4 \theta_3 (z_1-x_1)(z_1-y_1)(2z_1-x_1-y_1)$ and that
$\tC_0$ intersects $\Delta$ in  codimension one, at $\{x_1-y_1=0 \}$.

The  parametrization of $\tC_0$ in (\ref{eqn tC_0}) is   a  map
\bfo
\Upsilon:\R^5_{x_1, x_2, x_3, y_1, z_1}\times(\R_{\theta_3}\setminus 0)\lra \big(T^*\R^3\setminus0\big)\times \big(T^*\R^3\setminus0\big).
\efo
One easily sees
that $\Upsilon$ is singular at $\tS:=\Sigma_1(\Upsilon) =\{ x_1-y_1=x_1-z_1=0\}$, where  $\ker d\Upsilon=
\R\cdot\frac{\partial}{\partial z_1}$, which is
 $\not\subset T\Sigma_1$. Thus,  as discussed  below, $\tC_0$ is  an \emph{open umbrella}, exhibiting  the
 simplest kind of singularity of a Lagrangian manifold.

\subsection{Open umbrellas}\label{subsec ou}

For maps between manifolds of the same dimension, a   fold is a singularity of the type $S_{1,0}$  and a cusp is a
singularity of type $S_{1,1,0}$  \cite{Wh,Mo,GoGu}.  From larger dimensional spaces to smaller ones, $S_{1,0}$ maps are
submersions with folds,
while in the opposite direction they are referred to variously as \emph{cross caps} or  \emph{Whitney-Cayley umbrellas}  \cite{GoGu,Gi}.
In the lowest possible dimensions, if $g: R^2 \rightarrow R^3$ has a cross cap singularity then, in suitable local
coordinates, $(u,v,w)=g(x,y)=(x^2, y, xy)$, and its image is the algebraic surface $\{w^2=uv^2\}$,   the
Whitney-Cayley  umbrella. This is the simplest type of nonimmersed surface singularity in three
dimensions. It is an immersion away from the origin; it is actually an embedding on $\{y\ne0\}$; and along $\{y=0\}$, it is folded and hence 2-1.

Now, by adding one dimension to the range space, one may simultaneously both unfold the closed umbrella, making the
parametrization 1-1 away from
the origin, and  make it a Lagrangian manifold with singularity in $\R^4$. The \emph{open} (or unfolded) \emph{umbrella}  is the map $U
: R^2 \rightarrow R^4, \  U(x, y)=(x^2, y, xy, \frac{2}{3} x^3)$ (and its image). We have  $U^* \omega =0$, where $\omega$ is the
symplectic form on $R^4\simeq T^*\R^2$. Hence, the unfolded umbrella is a \emph{Lagrangian inclusion}, i.e., a variety which is a smoothly immersed Lagrangian manifold away from its singular points \cite{Gi}.  To find a similar conic model, one needs to
look  in $T^*R^3\backslash 0$;
for an unknown function $f(t,s)$,  consider
\[
\Lambda_1= \big\{ \big(s^2, t, f(t,s); st\theta, \frac{2}{3}t^3\theta, \theta\big)\in T^*\R^3\setminus 0:\theta\ne 0\, \big\}.
\]
Then, $\Lambda_1$ is a Lagrangian if $f$ is chosen so that  $d \xi \wedge d x|_{\Lambda_1}=0$. Since $d \xi \wedge d x= d(st\theta)
\wedge d(t^2) + d (\frac{2}{3} t^3\theta) \wedge ds + d \theta \wedge df= d\theta \wedge (2st^2 dt + \frac{2}{3} t^3 ds) + d
\theta \wedge df$, this holds if $df=- 2st^2 dt - \frac{2}{3} t^3 ds= - d(\frac{2}{3} t^3 s)$; using $f= -\frac{2}{3} t^3 s$
yields the conic Lagrangian,
$$\Lambda_1=\big \{\big( t^2, s, -\frac{2}{3}t^3s; st\theta, \frac{2}{3} t^3\theta, \theta\big):\, t,s\in\R,
\theta\in\R\backslash 0\, \big \},$$
exhibiting an open umbrella singularity along
$\{(0,0,0;0,0,\xi_3): \xi_3\ne 0\}$.
\medskip

\par Observe that  $\Lambda_1$ may be parametrized by a degenerate phase function. Letting $ \varphi (x, \theta_1, \theta_2, \tau, \eta)= (x_1-(\frac{\tau}{\eta})^2)
 \theta_1 + (x_3+\frac{2}{3} (\frac{\tau}{\eta})^3x_2)\theta_2 $,
 one has
 $$d_{\theta_1}\varphi= x_1-(\frac{\tau}{\eta})^2,\quad  d_{\theta_2}\varphi= x_3+ \frac{2}{3} (\frac{\tau}{\eta})^3 x_2  $$
 $$d_{\tau}\varphi=\frac{\tau}{\eta}(-2\frac{\theta_1}{\eta} + 2 \frac{\tau}{\eta} \frac{\theta_2}{\eta} x_2),\quad
 d_{\eta}\varphi= -(\frac{\tau}{\eta})^2(-2\frac{\theta_1}{\eta} + 2 \frac{\tau}{\eta} \frac{\theta_2}{\eta} x_2),   $$
 so that $\varphi$ is degenerate, with both $d_{\tau}\varphi$ and $d_{\eta}\varphi$ having normal crossings:
$d_{\tau}\varphi=0$ iff $\tau=0$ or $-2\frac{\theta_1}{\eta} + 2 \frac{\tau}{\eta} \frac{\theta_2}{\eta} x_2=0$, and
 $d_{\eta}\varphi=0$ on same two sets.
 For the second case, solving for $\theta_1= \frac{\tau}{\eta} x_2 \theta_2$ we see  that $\varphi$  parametrizes   $\Lambda_1$, while for the first case, $\varphi$  also parametrizes $\Lambda_0:=N^*\{ x_1=x_3=0 \}$. $\Lambda_0$ and $\Lambda_1$ intersect cleanly in codimension one, except at the singular set of $\Lambda_1$, which is contained in the intersection. $\varphi$  simultaneously parametrizes $\Lambda_0\cup\Lambda_1$, although it is not a multiphase function in the sense of \cite{Men}, used to parametrize a pair of cleanly intersecting smooth Lagrangians. The class of generalized FIOs we define in Sec. \ref{sec gfios} has similar features.
\medskip

\par For the general notion of an open umbrella, let $(M, \omega)$ be a smooth symplectic manifold of dimension $2n$.

\begin{definition}\label{def ou}
 $S \subset M$ is an \emph{open umbrella} if it is the image of a map $\psi: R^n \rightarrow M$ such that $d
 \psi$ drops rank by one simply at $\Sigma_1$, of codim two;  $\Ker (d \psi)|_{\Sigma_1} \not\subset T\Sigma_1\, $; and Range
 $(d\psi_p) $ is  Lagrangian for $p\in R^n \backslash \Sigma$.
\end{definition}

The class of open umbrellas is structurally stable under smooth perturbations in the class of Lagrangian inclusions \cite{Gi}.
As checked above, the parametrization $\Upsilon$ in (\ref{eqn tC_0}) satisfies the conditions in the definition, so that $\tilde{C}_0$ is
an open umbrella.

\bigskip

\subsection{Composition  for the model flat two sided cusp}\label{subsec comp}

\par We now want to show that the composition $A^*A$ for Fourier integral operators, such as $R_0$, associated with the model flat
two-sided cusp $C_0$ in (\ref{eqn czero})  results in operators with wave front relation in the union of the diagonal $\Delta$  and the
open umbrella $\tC_0$ from (\ref{eqn tC_0}). Let  $A\in I^m(\R^3,\R^3;C_0)$, so that $Af(x )=\int e^{i\phi(x, y, \theta)} a (x, y, \theta) f(y) d\theta  dy$ with $\phi(x,
y,\theta)=(x_2-y_2-(x_1-y_1)^2))\theta_2 \linebreak+ (x_3-y_3-(x_1-y_1)^4)\theta_3$ and $a \in S^{m+\frac12}_{1,0}$,  the standard H\"ormander class of symbols. We consider the normal operator formed by the composition,
$$A^*A f(x)=\int e^{i(\phi(z, y, \eta)-\phi(z, x, \theta))} a (z, y, \eta) \bar{a}(z, x, \theta) f(y) dy dz d \theta d \eta
.$$
We have
\bfo
\phi(z, y, \eta)-\phi(z, x, \theta) &=& (z_2-y_2-(z_1-y_1)^2))\eta_2 + (z_3-y_3-(z_1-y_1)^4)\eta_3 \\ & &-
(z_2-x_2-(z_1-x_1)^2))\theta_2 - (z_3-x_3-(z_1-x_1)^4)\theta_3.
\efo
After a stationary phase in $(z_2, \eta_2)$ and $(z_3, \eta_3)$  in the integral for $K_{A^*A}$, we obtain:
$$ \phi(z, y, \eta)-\phi(z, x, \theta) = (x_2-y_2+(z_1-x_1)^2-(z_1-y_1)^2)\theta_2 +
(x_3-y_3+(z_1-x_1)^4-(z_1-y_1)^4)\theta_3$$
$$= (x_2-y_2)\theta_2 + (x_3-y_3)\theta_3   + (y_1-x_1)[(2z_1-x_1-y_1)\theta_2
+(2z_1-y_1-x_1)((z_1-x_1)^2+(z_1-y_1)^2)\theta_3]$$
Let $z_1=\frac{\zeta}{\theta_3}$, and then consider the  change of variable,
$$\theta_1=-(2\frac{\zeta}{\theta_3}-x_1-y_1)\theta_2
-(2\frac{\zeta}{\theta_3}-y_1-x_1)((\frac{\zeta}{\theta_3}-x_1)^2+(\frac{\zeta}{\theta_3}-y_1)^2)\theta_3,   $$
for which
\bfo
\frac{\partial \theta_1}{\partial \zeta}&=& -2 \frac{\theta_2}{\theta_3}- 2
((\frac{\zeta}{\theta_3}-x_1)^2+(\frac{\zeta}{\theta_3}-y_1)^2)) - 2 (2\frac{\zeta}{\theta_3} -x_1-y_1)^2,\\
\frac{\partial^2 \theta_1}{\partial \zeta^2}&=& -\frac{12}{\theta_3}(2\frac{\zeta}{\theta_3}-x_1-y_1), \hbox{ and }\quad
\frac{\partial^3 \theta_1}{\partial \zeta^3} =- \frac{24}{\theta_3}.
\efo
When $\frac{\partial^2 \theta_1}{\partial \zeta^2}=0$ then $2\frac{\zeta}{\theta_3}=x_1+y_1$,  $\theta_1=0$ and
$\theta_2=-\frac{(x_1-y_1)^2}{2}\theta_3$.
  Let
  $T(t_1, t_2)=(t_1 t_2+\frac{1}{2}t_1^3, t_2)$ be a normal form for  a two-dimensional map with a cusp singularity. To
  understand the push forward $T_*u$  of a distribution $u$ under the  map $T$, note that the kernel of $T_{*}$ is
  $$K_{T_{*}}(z,t)=\delta(z-T(t)) =\int e^{i[(z_1-t_1t_2-\frac12t_1^3)\sigma_1+(z_2-t_2)\sigma_2]}  1 (\sigma) d\sigma, $$
  so that
   $$T_{*}u(z)= \int e^{i [z_1\sigma_1+z_2\sigma_2]}\Big [ \int e^{-i[t_2\sigma_2+(t_1t_2+\frac12t_1^3)\sigma_1]} u(t)d t\Big ]d \sigma.$$
 \par Returning now to $A^*A$, letting $z_1=\frac{\theta_1}{\theta_3}$ and $z_2=\frac{\theta_2}{\theta_3}+\frac{1}{2}(x_1-y_1)^2$,
 the phase function of $A^*A$ becomes
 \bfo
 \tilde{\phi}(x, y, \theta, \sigma)&\!\!\!=\!\!\!&(x_1-y_1)\theta_1+(x_2-y_2)\theta_2+(x_3-y_3)\theta_3 \\
 & &\qquad\qquad\quad +(\frac{\theta_1}{\theta_3}-t_1t_2-\frac{1}{2}t_1^3) \sigma_1 +
 (\frac{\theta_2}{\theta_3}+\frac{1}{2}(x_1-y_1)^2-t_2)\sigma_2
 \efo
 and the amplitude becomes: $a\times 1(\sigma)$.
 Next, we perform stationary phase, first in $t_2, \sigma_2$  and then in $ \theta_1, \sigma_1$. We have
 \bfo
 d_{t_2} \tilde{\phi} = -t_1\sigma_1-\sigma_2,\quad
 d_{\sigma_2}\tilde{\phi} = \frac{\theta_2}{\theta_3} + \frac{1}{2} (x_1-y_1)^2-t_2.
 \efo
 The Hessian is $1$, the  phase function is $\tilde{\phi} = (x-y)\theta +(
 \frac{\theta_1}{\theta_3}-\frac{1}{2}t_1^3-t_1[\frac{\theta_2}{\theta_3} + \frac{1}{2}(x_1-y_1)^2])\sigma_1$, and
 \bfo
 d_{\theta_1}\tilde{\phi} = x_1-y_1+ \frac{\sigma_1}{\theta_3}, \quad
 d_{\sigma_1}\tilde{\phi} = \frac{\theta_1}{\theta_3} -\frac{1}{2} t_1^3-t_1(\frac{\theta_2}{\theta_3} +
 \frac{1}{2}(x_1-y_1)^2).
 \efo
 The Hessian is $\frac{1}{\theta_3^2}$. Finally, we introduce the variable $\tau= t_1 \theta_3$  and  obtain
 \beq\label{eqn modeltp}\nonumber
   \tilde{\phi} = \left(x_2-y_2+\frac{\tau}{\theta_3}(x_1-y_1)\right)\theta_2
   +\left(x_3-y_3+\frac{1}{2}(x_1-y_1)(\frac{\tau}{\theta_3})^3+ \frac{1}{2}
   \frac{\tau}{\theta_3}(x_1-y_1)^3\right)\theta_3,
   \eeq
with respect to which   the Schwartz
kernel has the representation
$$K_{A^*A}(x,y)=\int e^{i \tilde{\phi}} b\, d\theta_2 d\theta_3 d\tau,\quad b\in S^{2m}.$$
The wave front relation satisfies
 \bfo
WF(K_{A^*A})'&\subset&\Big\{ \big(x_1, x_2, x_3, \frac{\tau\theta_2}{\theta_3} +\frac{1}{2}\frac{\tau^3}{\theta_3^2} +
\frac{3}{2} (x_1-y_1)^2 \tau, \theta_2, \theta_3;  \\
 & &\qquad y_1, y_2, y_3, \frac{\tau\theta_2}{\theta_3} +\frac{1}{2}\frac{\tau^3}{\theta_3^2}+ \frac{3}{2} (x_1-y_1)^2 \tau,
 \theta_2, \theta_3\big): \\
& &\qquad d_{\theta_2} \tilde{\phi}= x_2-y_2+\frac{\tau}{\theta_3}(x_1-y_1)=0, \\
& &\qquad  d_{\theta_3} \tilde{\phi}=x_3-y_3-\frac{\theta_2\tau}{\theta_3^2}(x_1-y_1)-\frac{\tau^3}{\theta_3^3}(x_1-y_1) =0,
\\
& &\qquad   d_{\tau} \tilde{\phi}= \frac{\theta_2}{\theta_3}(x_1-y_1) + \frac{3}{2} \frac{\tau^2}{\theta_3^2}(x_1-y_1)+
\frac{1}{2} (x_1-y_1)^3=0\Big\}.
\efo
From the last relation, we see that the phase function is degenerate, with the critical set in the phase variables  having a
normal crossing,  being the union of two transverse surfaces, corresponding to  $x_1-y_1=0$ or $\frac{\theta_2}{\theta_3}+
\frac{3}{2} \frac{\tau^2}{\theta_3^2} + \frac{1}{2} (x_1-y_1)^2=\nolinebreak0$. The points where $x_1=y_1$ give rise to the
diagonal $\Delta$, but with the parametrization actually being a cusp:
$$(x,\tau,\theta)\lra \big(x, \frac{\tau\theta_2}{\theta_3} +\frac{1}{2} \frac{\tau^3}{\theta_3^2}, \theta_2, \theta_3; \ \
x,  \frac{\tau\theta_2}{\theta_3} +\frac{1}{2} \frac{\tau^3}{\theta_3^2}, \theta_2, \theta_3   \big). $$
On the other hand, the points where $\frac{\theta_2}{\theta_3}+ \frac{3}{2} \frac{\tau^2}{\theta_3^2} +
\frac{1}{2}(x_1-y_1)^2=0$ contribute to $WF(A^*A)$ the set, parametrized by $x\in\R^3, y_1,\tau\in\R, \tht\ne 0$,
\beq\label{eqn tCo}
{\tC_{0}}=\Big\{ \quad\big(x, \xi;\,  y, \eta\big): \,  y_2&=& x_2+\frac{\tau}{\theta_3}(x_1-y_1),\quad \xi_3=\eta_3=\theta_3, \\
y_3&=&  x_3+\frac{1}{2}\frac{\tau}{\theta_3}(x_1-y_1)^3 +\frac{1}{2} \frac{\tau^3}{\theta_3^3} (x_1-y_1),\nonumber \\
\xi_1&=&\eta_1= -\frac{\tau^3}{\theta_3^2}+ (x_1-y_1)^2\tau, \nonumber  \\
\xi_2&=&\eta_2= -\frac{1}{2}(x_1-y_1)^2\theta_3 -\frac{3}{2} \frac{\tau^2}{\theta_3}\Big\}, \nonumber
\eeq
which is the image of $\rho(x_1, x_2, x_3, y_1, \theta_3, \tau)=(x, \xi; \ \ y, \eta) $ with $y_2, y_3, \eta_1, \eta_2,
\eta_3$, $\xi_1, \xi_2, \xi_3$ given by the above relations. $\rho$ satisfies the conditions of  Def. \ref{def ou}:  $d\rho$ drops rank simply by one at  $\Sigma=\{x_1-y_1=\tau=0\}$ and  $\Ker d\rho=\R\cdot\frac{\partial}{\partial \tau} \nsubseteq T\Sigma$.
Hence $\tC_{0}$ is an open umbrella.

\par Next, we  compute the principal symbols on $\Delta$ and $\tC_0$ away from their intersection.
We have $Crit_{\tilde{\phi}}= Crit_{\Delta} \cup Crit_{\tC_0}= \{ (x, y, \tau, \theta_2, \theta_3)\,| d_{\tau} \tilde{\phi}=0,
d_{\theta_2} \tilde{\phi}=0, d_{\theta_3} \tilde{\phi}=0 \} $. Using  H\"{o}rmander's formula \cite{H-book}, $\sigma=a
(E_{\tilde{\phi}})^{\frac{1}{2}}$ where $E_{\tilde{\phi}}=|\frac{D(\lambda_i, \frac{\partial \tilde{\phi}}{\partial
\theta})}{D(x, \theta)}|^{-1}$ and $\lambda_i$ are  local coordinates  on $Crit_{\tilde{\phi}}$. On $Crit_{\Delta}$,
local coordinates are $(x,\tau, \theta_2, \theta_3)$ and $E_{\tilde{\phi}}=|\frac{D(x,\tau, \theta_2, \theta_3 ,
\frac{\partial \tilde{\phi}}{\partial \tau}, \frac{\partial \tilde{\phi}}{\partial \theta_2}, \frac{\partial
\tilde{\phi}}{\partial \theta_3})}{D(x, \tau, \theta_2, \theta_3, y_1, y_2, y_3)}|^{-1}=( \frac{\theta_2}{\theta_3} +
\frac{3}{2}\frac{\tau^2}{\theta_3^2})^{-1}$, while, on $Crit_{\tC_0}$,  local coordinates are $(x,y_1, \tau, \theta_3)$ and
$E_{\tilde{\phi}}=|\frac{D(x, y_1, \tau,  \theta_3 , \frac{\partial \tilde{\phi}}{\partial \tau}, \frac{\partial
\tilde{\phi}}{\partial \theta_2}, \frac{\partial \tilde{\phi}}{\partial \theta_3})}{D(x, y_1, \tau, \theta_2, \theta_3, y_2,
y_3, \theta_2)}|^{-1}=( \frac{x_1-y_1}{\theta_3})^{-1}$. Thus, both principal symbols
are singular as one approaches $\Delta\cap \tC_{0}$, behaving as $\delta^{-\frac12}$,
where $\delta$ is the distance on either $\Delta$ or $\tC_{0}$  to the intersection.

\subsection{Flat two-sided cusps}\label{subsec class}

Motivated by the previous example, we now define the class of canonical relations for which we will establish a composition
calculus.

\begin{definition}\label{def ftsc}
If $X$ and $Y$ are  manifolds of dimension $n\ge 3$, then a  canonical relation  $C\subset (T^*X\setminus 0)\times
(T^*Y\setminus 0)$  is a \emph{flat two-sided cusp} if
\medskip

\ni (i) both $\pl:C\lra T^*X$ and $\pr:C\lra T^*Y$ have at most cusp singularities;

\ni (ii) the   left- and right-cusp points are equal: $\Soo(\pl)=\Soo(\pr)=:\Soo$; and

\ni (iii) $\pl(\Soo)\subset T^*X$ and $\pr(\Soo)\subset T^*Y$ are coisotropic (involutive) and nonradial.
\end{definition}

Here, as usual, \emph{at most cusp} means either   cusps or folds, which are unavoidable in the neighborhood of a cusp.
\emph{Nonradial}, meaning that the the radial vector field $\sum \xi_j \partial_{\xi_j}$ does not lie in  the two-dimensional
Hamiltonian foliation of the codimension two coisotropic submanifold, is a  standard technical assumption   needed to apply the homogeneous
Darboux theorem in Sec. \nolinebreak\ref{sec cusps}.
\bigskip

Other examples of  flat two-sided cusps coming from  generalized Radon transforms, verified as with the example above, are
those for translates of  curves  in $\R^4$ satisfying (\ref{eqn rfour}), and  for the translates in $\R^3$ of  Menn's surface \cite{Ban},
$$C_1=N^{*} \{ x_3-y_3=(x_1-y_1)^2(x_2-y_2) - (x_2-y_2)^2  \}' , $$
or any other surface in $\R^3$  whose Gauss map has a cusp.
Computation also shows that the composition  $C^t_1 \circ C_1 = \Delta \cup \tC_1$, with $\tC_1$ similar to $\tC_0$ as in (\ref{eqn tCo}) above.
\medskip

\noindent{\bf Remarks. 1.}  If (i) and (ii) in Def. \ref{def ftsc} are satisfied, the structure of
$\omega_{T^*X}|_{T\pl(\Soo)}$ is the same as $\omega_{T^*Y}|_{T\pr(\Soo)}$, since both are equivalent to the restriction to
$\Soo\subset C$ of the degenerate symplectic form $\omega_C:=\pl^*\omega_{T^*X}=\pr^*\omega_{T^*Y}$. Thus, $\pl(\Soo)$ is involutive
iff $\pr(\Soo)$ is, and one can see that this holds iff the two-dimensional $\ker(\omega_C|_{\So})$ is simply tangent to $\So$ at $\Soo$. This is in contrast to the situation when the image of the cusp set is symplectic (strongly noninvolutive), in which case $\omega_C$ is a folded symplectic form and $\ker(\omega_C|_{\So})\pitchfork\So$ everywhere \cite{M}.

\noindent{\bf 2.} Condition (i) is stable under small perturbations, since cusps are structurally stable \cite{GoGu}.
However, given (i), condition (ii)  is unstable and atypical. To see this concretely, consider (possibly) nontranslationally
invariant families of curves $\{\gamma_x\}_{x\in\R^4}$, using the framework of \cite{CNSW}  as analyzed in terms of FIOs with cusps
in \cite{GS}. For  vector fields $X,Y,Z,W$ on $\R^4$, with $X\ne 0$, let
\bfo
\gamma_x(t)=exp_x(tX+t^2Y+t^3Z+t^4W).
\efo
The associated generalized Radon transform is in $I^{-\frac12}(\R^4,\R^4;C)$,
with $C$  a canonical relation  for which, by \cite[Prop. 6.1]{GS},  $\pl$ (resp., $\pr$) is at most a cusp if
\bfo
X,\, Y,\, Z\pm\frac16[X,Y],\, W\pm\frac14[X,Z]+\frac1{24}[X,[X,Y]]\hbox{ are linearly independent.}
\efo
From the calculations in \cite{GS}, one sees that  a necessary condition for $\Soo(\pl)=\Soo(\pr)$ is that $(Z +\frac16[X,Y])-(Z-\frac16[X,Y])\in\hbox{
span }\{X,Y\}$, i.e.,  $\hbox{ span }\{X,Y\}$ is an integrable distribution of 2-planes, which  is generically
not the case.

Nevertheless, as seen in the next section,  flat two-sided cusps
arise naturally in a  quite different setting without translation invariance or integrability.

\section{Linearized seismic inversion with cusp caustics}\label{sec caustics}

\par We consider linearized seismic inversion for the  single source data set,  under the assumption that  the background
sound speed has  caustics of at most cusp type, i.e., either folds or cusps. We show that the linearized forward scattering
operator $F$ is a Fourier integral operator associated with a flat two-sided cusp, $\tC$. We  only briefly describe this
problem, referring to \cite{NS-fold,NS,N,F,FG} for more motivation and details.

Acoustic waves are generated at the surface of the earth, scatter off heterogeneities in the
subsurface and return to the surface, where measurements of the pressure field  are used to reconstruct an
image of the subsurface.
The model for the scattered waves is given by the acoustic wave equation,
\beq\label{eqn awe}
 \frac{1}{c^2(x)} \frac{\partial^2p}{\partial
t^2}(x,t)-\triangle p (x,t)=\delta(x-s)\, \delta(t),\qquad
 p(x,t)=0,\quad  t < 0,
\eeq  where $x \in Y= \{x \in R^3, x_3 \geq 0 \}$, representing the Earth, $p(x,t)$
is the pressure field resulting from a pulse at  the source $s$ at time $t=0$,
and $c(x)$ is the unknown sound speed field.

To make this nonlinear inverse problem tractable, one considers a linearized
operator $F$ which maps singular perturbations of a smooth
background  sound speed in the subsurface, assumed known, to
perturbations of the resulting pressure field at the surface.
Thus, the linearization consists in
assuming $c$ to be of the form $c=c_0 + \delta c$ and the resulting $p=p_0+\delta p$, where $c_0 $ is
a smooth known background field. The formal linearization of (\ref{eqn awe}) is
\beq\label{eqn lin}
\quad\square_{c_0}(\delta p):= \frac{1}{c_0^2(x)}\frac{\partial^2 \delta
p}{\partial  t^2}(x,t) -\triangle \delta p(x,t)= \frac{2}{c_0^3 } \frac{\partial^2 p_0}{\partial t^2}\cdot\delta
c(x),\quad
\delta p =0, \quad t < 0,
\eeq
where
$p_0$, the Green's function for $\square_{c_0}$.
The \emph{linearized
scattering operator}  is   $F : \delta c \rightarrow \delta
p|_{\partial Y \times (0,T)}$.  Under mild technical assumptions, $F$ is an FIO associated with a canonical relation $C\subset (T^*X\setminus 0)\times (T^*Y\setminus 0)$, where $X=\partial Y\times (0,T)$ is the data space \cite{R,tKSV,NS}.
The
goal is then to left-invert $F$; standard techniques suggest  studying
left invertibility of the normal operator, $N=F^*F$.

\par Let $H(x,\xi)=\frac{1}{2}(c_0(x)^{-2}-|\xi|^2)$  be the Hamiltonian
associated to  $c_0$, and $\L_s$ the image of
$T_s^*\R^3
\setminus 0$ under the bicharacteristic flow  associated to $H$,
which is a Lagrangian submanifold of $T^*\R^3 \setminus 0$.
A \emph{caustic}  is a singularity of
the spatial projection $\pi : \L_s \rightarrow Y$.
It is known that the  caustics exhibited for generic soundspeeds  $c_0$ are the same as generic Lagrangian singularities, i.e., folds, cusps, swallowtails, etc. \cite{Guck}.

If $c_0$ is such that only fold caustics occur,
it was shown by Nolan \cite{N} that $F\in I^1(X,Y;C)$ with $C$ a two-sided fold.  Furthermore, the
Schwartz kernel of the operator $F^*F$ belongs to a class of
distributions associated to two cleanly intersecting Lagrangians
in $(T^*Y \setminus 0) \times (T^*Y \setminus 0)$, with the corresponding canonical relations being the
diagonal $\Delta$ and another folding canonical relation \cite{N,F}.

The next caustics to consider are those of cusp type,
 meaning that the only singularities of
the spatial projection $\pi : \L_s \rightarrow Y$ are folds or cusps.
A cusp caustic is already present in the 2D example of \cite{NS-fold} exhibiting loss of regularity for $F$.
To analyze the geometry of the canonical relation $C$ in the presence of cusp caustics, we make
use of the description of $\L_s$ in \cite{N}. It can be
parametrized by $t_{inc}$, the time travelled by the incident
ray, and the takeoff direction $(p_1,p_2,p_3) \in \mathbb S ^2$. We
can change these coordinates to $(x_1, x_2, p_3)$. Hence on
$\tC_s$, $x_3=f(x_1, x_2, p_3)$ and $(p_1,p_2)=(g_1(x_1, x_2,
p_3), g_2(x_1, x_2, p_3))$ and $\L_s$ is the graph of a function $\nabla G(x_1, x_2, p_3)$ which means that $\frac{\partial
G}{\partial x_1}=g_1, \frac{\partial G}{\partial x_2}=g_2,\frac{\partial G}{\partial p_3}=f$. Then $\frac{\partial
g_1}{\partial p_3}=\frac{\partial f}{\partial x_1}, \frac{\partial g_2}{\partial p_3}=\frac{\partial f}{\partial x_2},
\frac{\partial g_2}{\partial x_1}=\frac{\partial g_1}{\partial x_2}$. In this new setting, $\pi(x_1,x_2,p_3)=(x_1, x_2,
f(x_1,x_2,p_3))$
det $d \pi=\frac{\partial f}{\partial p_3}$ and  cusp
caustics occur when
\beq\label{eqn caustic}
\quad\frac{\partial f}{\partial p_3} =\frac{\partial^2 f}{\partial p_3^2}= 0, \quad \frac{\partial^3 f}{\partial p_3^3}
\neq 0, \hbox{  and }
\big\{\nabla \frac{\partial f}{\partial p_3},\, \nabla \frac{\partial^2 f}{\partial p_3^2}\big\}\hbox{ are linearly
independent.}
\eeq

Next, we  parametrize the canonical relation $C$ of $F$ in
terms of $ x_1, x_2$ and
$p_3$ ; $(\alpha_1, \alpha_2, \sqrt {1-|\alpha|^2})$, the take off
direction of the reflected ray, writing $\alpha=
(\alpha_1,\alpha_2)$; and $\tau$, the variable dual  to time.
\bfo
C& &= \big\{ ( r_1(\cdot), r_2(\cdot),
t_{inc}(\cdot) + t_{ref}(\cdot),
\rho_1(\cdot),   \rho_2(\cdot), \tau ; \, x_1, x_2, f(\cdot),\\
& & \qquad-\tau(c_0^{-1}(\cdot)
\alpha_1 + g_1(\cdot)),
 -\tau (c_0^{-1}(\cdot)\alpha_2 + g_2(\cdot)),-\tau(c_0^{-1}(\cdot) \sqrt{1-|\alpha|^2} + p_3) ) \big\},
\efo
where
\bfo
f(\cdot)&=&f(x_1,x_2, p_3); \qquad
r_j(\cdot)=r_j(x_1,x_2, f(x_1,x_2, p_3),\alpha),  j=1,2;\\
t_{inc}(\cdot)&=&t_{inc}(x_1, x_2, p_3);\qquad
t_{ref}(\cdot)=t_{ref}(x_1,x_2, f(x_1,x_2,p_3), \alpha);\\
\rho_j(\cdot)&=&\rho_j(x_1,x_2, f(x_1,x_2,  p_3), \alpha),  j=1,2;\qquad
g_j(\cdot)=g_j(x_1, x_2, p_3),  j=1,2; \\
& & \hbox{ and }  c_0^{-1}(\cdot)=c_0^{-1}(x_1, x_2,f(\cdot)).
\efo
We have $\pi_R(x_1, x_2, p_3, \alpha_1, \alpha_2, \tau)=(x_1, x_2, f(\cdot), -\tau(c_0^{-1}(\cdot)
\alpha_1 + g_1(\cdot)),  -\tau (c_0^{-1}(\cdot)\alpha_2 + g_2(\cdot)),
-\tau(c_0^{-1}(\cdot)\sqrt{1-|\alpha|^2} + p_3)$
and thus
$$d\pi_R= \left[ \begin{array} {cccccc} 1 & 0 & 0 & 0 & 0 & 0  \\
0 &  1  & 0 & 0 & 0 & 0  \\
\frac{\partial f}{\partial x_1} & \frac{\partial f}{\partial x_2} & \frac{\partial f}{\partial p_3} & 0 & 0 & 0  \\
A & B & C & -\tau c_0^{-1}  & 0 &  -c_0^{-1}(\cdot)\alpha_1 - g_1(\cdot)\\
D & E & F & 0 & -\tau c_0^{-1} &  -c_0^{-1}(\cdot)\alpha_2 -g_2(\cdot)\\
G & H & I & \frac{\tau c_0^{-1} \alpha_1}{\sqrt{1-|\alpha|^2}} & \frac{\tau c_0^{-1} \alpha_2}{\sqrt{1-|\alpha|^2}} &
-c_0^{-1}(\cdot)\sqrt{1-|\alpha|^2} - p_3 \\
\end{array} \right]$$
where $A=\frac{\partial(-\tau(c_0^{-1}(\cdot)
\alpha_1 + g_1(\cdot)))}{\partial x_1}$; $B=\frac{\partial(-\tau(c_0^{-1}(\cdot)
\alpha_1 + g_1(\cdot)))}{\partial x_2}$; $C=\frac{\partial(-\tau(c_0^{-1}(\cdot)
\alpha_1 + g_1(\cdot)))}{\partial p_3}$; $D=\frac{\partial(-\tau (c_0^{-1}(\cdot)\alpha_2 + g_2(\cdot)))}{\partial x_1}$;
$E=\frac{\partial(-\tau (c_0^{-1}(\cdot)\alpha_2 + g_2(\cdot)))}{\partial x_2}$; $F=\frac{\partial(-\tau
(c_0^{-1}(\cdot)\alpha_2 + g_2(\cdot)))}{\partial p_3}$; $G=\frac{\partial (-\tau(c_0^{-1}(\cdot)\sqrt{1-|\alpha|^2} +
p_3))}{\partial x_1}$; $H=\frac{\partial (-\tau(c_0^{-1}(\cdot)\sqrt{1-|\alpha|^2} + p_3))}{\partial x_2}$; and
$I=\frac{\partial (-\tau(c_0^{-1}(\cdot)\sqrt{1-|\alpha|^2} + p_3))}{\partial p_3}$.
Hence  $\det\, d \pi_R= \frac{\partial f}{\partial p_3}$ and the kernel at $\So=\{\det\, d \pi_R= 0\}$ is
spanned by  a vector $v_1$ which is $\frac{\partial}{\partial p_3}$ plus
 a linear combination of $\{\frac{\partial}{\partial \alpha_1}, \frac{\partial}{\partial \alpha_2},\frac{\partial}{\partial \tau} \}$, because    the matrix
$$\left[ \begin{array} {ccc}
-\tau c_0^{-1}  & 0 &  c_0^{-1}(\cdot)\alpha_1 + g_1(\cdot)\\
  0 & -\tau c_0^{-1} &  c_0^{-1}(\cdot)\alpha_2 + g_2(\cdot)\\
 \frac{\tau c_0^{-1} \alpha_1}{\sqrt{1-|\alpha|^2}} & \frac{\tau c_0^{-1} \alpha_2}{\sqrt{1-|\alpha|^2}} &
 c_0^{-1}(\cdot)\sqrt{1-|\alpha|^2} + p_3 \\
\end{array} \right]$$  is nondegenerate by the calculations of  \cite{N}.
From the cusp conditions (\ref{eqn caustic}), one can  see that $\pi_R$ has a cusp singularity.
\par Similarly,  $\pi_L(x_1, x_2, \alpha_1,\alpha_2,p_3,\tau)=(r_1, r_2, t_{inc}+t_{ref}, \rho_1, \rho_2, \tau)$ and
$$d\pi_L= \left[ \begin{array} {cccccc} \frac{\partial r_1}{\partial x_1} + \frac{\partial r_1}{\partial x_3} \frac{\partial
f}{\partial x_1} & \frac{\partial r_1}{\partial x_2} + \frac{\partial r_1}{\partial x_3} \frac{\partial f}{\partial x_2}  &
\frac{\partial r_1}{\partial \alpha_1} & \frac{\partial r_1}{\partial \alpha_2} & \frac{\partial r_1}{\partial x_3}\frac{
\partial f}{\partial p_3} & 0  \\
 \frac{\partial r_2}{\partial x_1} + \frac{\partial r_2}{\partial x_3} \frac{\partial f}{\partial x_1} & \frac{\partial
 r_2}{\partial x_2} + \frac{\partial r_2}{\partial x_3} \frac{\partial f}{\partial x_2}  & \frac{\partial r_2}{\partial
 \alpha_1} & \frac{\partial r_2}{\partial \alpha_2} & \frac{\partial r_2}{\partial x_3}\frac{ \partial f}{\partial p_3} & 0
 \\
 \frac{\partial t_{ref}}{\partial x_1} + \frac{\partial t_{inc}}{\partial x_1} & \frac{\partial t_{ref}}{\partial x_2} +
 \frac{\partial t_{inc}}{\partial x_2} & \frac{\partial t_{ref}}{\partial \alpha_1} & \frac{\partial t_{ref}}{\partial
 \alpha_2} & \frac{\partial t_{ref}}{\partial x_3}\frac{ \partial f}{\partial p_3} + \frac{\partial t_{inc}}{\partial p_3}& 0
 \\
\frac{\partial \rho_1}{\partial x_1} + \frac{\partial \rho_1}{\partial x_3} \frac{\partial f}{\partial x_1} & \frac{\partial
\rho_1}{\partial x_2} + \frac{\partial \rho_1}{\partial x_3} \frac{\partial f}{\partial x_2}  & \frac{\partial
\rho_1}{\partial \alpha_1} & \frac{\partial \rho_1}{\partial \alpha_2} & \frac{\partial \rho_1}{\partial x_3}\frac{ \partial
f}{\partial p_3} & 0  \\
\frac{\partial \rho_2}{\partial x_1} + \frac{\partial \rho_2}{\partial x_3} \frac{\partial f}{\partial x_1} & \frac{\partial
\rho_2}{\partial x_2} + \frac{\partial \rho_2}{\partial x_3} \frac{\partial f}{\partial x_2}  & \frac{\partial
\rho_2}{\partial \alpha_1} & \frac{\partial \rho_2}{\partial \alpha_2} & \frac{\partial \rho_2}{\partial x_3}\frac{ \partial
f}{\partial p_3} & 0  \\
0 & 0 & 0 & 0 & 0 & 1 \\
\end{array} \right]$$
From \cite{N} we have that $\frac{\partial t_{inc}}{\partial p_3}=0, \frac{\partial r_1}{\partial x_3}=0, \
\frac{\partial r_2}{\partial x_3}=0, \ \frac{\partial \rho_1}{\partial x_3}=0, \ \frac{\partial \rho_2}{\partial x_3}=0 , \
\frac{\partial t_{ref}}{\partial x_3}<0 $,
and that the matrix $$\left[ \begin{array} {cccc}
 \frac{\partial r_1}{\partial x_1} & \frac{\partial r_1}{\partial x_2} &  \frac{\partial r_1}{\partial \alpha_1} &
 \frac{\partial r_1}{\partial \alpha_2}\\
  \frac{\partial r_2}{\partial x_1} & \frac{\partial r_2}{\partial x_2} &  \frac{\partial r_2}{\partial \alpha_1} &
  \frac{\partial r_2}{\partial \alpha_2}\\
  \frac{\partial \rho_1}{\partial x_1} & \frac{\partial \rho_1}{\partial x_2} &  \frac{\partial \rho_1}{\partial \alpha_1} &
  \frac{\partial \rho_1}{\partial \alpha_2} \\
  \frac{\partial \rho_2}{\partial x_1} & \frac{\partial \rho_2}{\partial x_2} &  \frac{\partial \rho_2}{\partial \alpha_1} &
  \frac{\partial \rho_2}{\partial \alpha_2} \\
\end{array} \right]$$  is nondegenerate.
Thus,  $ \det\, d \pi _L=\frac{\partial f}{\partial p_3}$,  $\ker\pi_L=\frac{\partial}{\partial p_3}$ and as before one sees
that $\pi_L$ has also a cusp singularity.

\par Now we will show that the images of the cusp points are involutive. By Remark 1 in Sec. \ref{subsec class},  it is enough to examine the
degenerate symplectic form $\omega_C=\pi_R^* \omega_{T^*Y}$  on $C$
and its kernel, equal to $\{v: \omega_C(v,w)=0, \forall w \in
TC  \}$.  We have from  (\ref{eqn caustic}) that $\Sigma_1=\{ (x,\xi);  f_{p_3}=0 \}$ and
$\Sigma_{1,1}=\{ (x,\xi); f_{p_3}=f_{p_3p_3}=0 \}$.
Since $\omega_C|_{\Sigma_1}$ drops rank by $2$ and  $\ker d\pi_R= \R\cdot v_1$, one has
$\ker\omega_C=\hbox{ span}\{  v_1, v_2 \}$, where $v_2$ is a  vector, without $\frac{\partial}{\partial p_3}$ component, that we now determine.
At $\So$, one has
\bfo
\omega_C|_{\Sigma_1}&=& (c_0^{-1} \alpha_1 + g_1 + (c_0^{-1}\sqrt {1-|\alpha|^2 }+p_3)f_{x_1})  (d\tau \wedge d x_1)
 \\
&+&  (c_0^{-1} \alpha_2 + g_2 + (c_0^{-1}\sqrt {1-|\alpha|^2 } +p_3)
f_{x_2}) (d\tau \wedge d x_2)\\
&+& \tau c_0^{-1}(1-\frac{\alpha_1}{\sqrt{1-|\alpha|^2}}
f_{x_1})(d \alpha_1 \wedge dx_1)+  \tau c_0^{-1}(1-\frac{\alpha_2}{\sqrt{1-|\alpha|^2}} f_{x_2})(d \alpha_2 \wedge dx_2)  \\
&+&  \tau c_0^{-1} \frac{-\alpha_1}{\sqrt{1-|\alpha|^2}} f_{x_2}(d\alpha_1 \wedge dx_2)+  \tau c_0^{-1} \frac{-\alpha_2}{\sqrt{1-|\alpha|^2}} f_{x_1}(d\alpha_2 \wedge dx_1) \\
&+&  \tau (\partial_{x_2}c_0^{-1} \alpha_1-\partial_{x_1}c_0^{-1} \alpha_2 + \partial_{x_2}
c_0^{-1}\sqrt{1-|\alpha|^2} f_{x_1} -\partial_{x_1} c_0^{-1}\sqrt{1-|\alpha|^2} f_{x_2})(dx_2 \wedge dx_1) \\
&+&  2 \tau f_{x_1}(d p_3 \wedge dx_1) +  2 \tau f_{x_2}(d p_3 \wedge dx_2).
\efo
We need to look at the $5 \times 5$ skew symmetric matrix of
$\omega_C$ corresponding to $x_1,\! x_2,\! \alpha_1,\! \alpha_2,\! \tau\!:$
$$ \left[ \begin{array} {ccccc} 0 & a & \star &  \star & -b  \\
-a & 0 & \star & \star &  -c  \\
 \tau c_0^{-1}(1-\frac{\alpha_1}{\sqrt{1-|\alpha|^2}} f_{x_1})& \tau c_0^{-1}
 \frac{-\alpha_1}{\sqrt{1-|\alpha|^2}} f_{x_2} & 0 & 0 & 0  \\
 \tau c_0^{-1} \frac{-\alpha_2}{\sqrt{1-|\alpha|^2}} f_{x_1}& \tau
 c_0^{-1}(1-\frac{\alpha_2}{\sqrt{1-|\alpha|^2}} f_{x_2}) &   0  & 0 & 0 \\
b & c & 0 & 0 & 0 \\
\end{array} \right]$$
where $a$ is the coefficient of $d x_1 \wedge dx_2$ , $b$ is the coefficient of  $d \tau \wedge d x_1$ and  $c$ is the
coefficient of $ d \tau \wedge dx_2$. $v$ is a combination of $\frac{\partial }{\partial x_1}, \frac{\partial
}{\partial x_2},\frac{\partial }{\partial \alpha_1},\frac{\partial }{\partial \alpha_2},\frac{\partial }{\partial \tau}$;
note that
\bfo
 \left| \begin{array} {cc}
1-\frac{\alpha_1}{\sqrt{1-|\alpha|^2}} f_{x_1}&
 \frac{-\alpha_1}{\sqrt{1-|\alpha|^2}} f_{x_2} \\
\frac{-\alpha_2}{\sqrt{1-|\alpha|^2}} f_{x_1}&
1-\frac{\alpha_2}{\sqrt{1-|\alpha|^2}} f_{x_2}  \\
\end{array} \right|
 =1+ \nabla {\sqrt{1-|\alpha|^2}}  \cdot \nabla_{x_1,x_2} f \neq 0,
\efo
which means that $v_2$ does not contain any $\frac{\partial }{\partial x_1}, \frac{\partial }{\partial x_2}$ terms, but only
$\frac{\partial }{\partial \alpha_1},\frac{\partial }{\partial \alpha_2},\frac{\partial }{\partial \tau}$. This shows that
$\ker\omega_C$ is simply tangent to $\Sigma_1$ at  $\Sigma_{1,1}$, and thus $\pr(\Soo)$ is involutive; by Remark 1, so is $\pl(\Soo)$.

\section{Weak normal forms for flat two-sided cusps }\label{sec cusps}

\par Next, we show that any flat two-sided cusp can be prepared, by application of suitable canonical transformations on the
left and right, so as to be parametrized by a phase function similar to that for $R_0$  in (\ref{eqn Rnot}).  Recall that one can conjugate any two-sided \emph{fold} to a single normal form \cite{M,meta}, but, as in \cite{GU1,FG}, we will need to work with merely approximate normal forms. For simplicity, we restrict ourselves to the three dimensional setting of interest for the seismic problem.

Assume $\dim(X)=\dim(Y)=3$ and suppose  $C\subset T^*X\times T^*Y$ is a flat two-sided cusp.
 Let $c^0\in\Sigma_{1,1}$. We claim that there exist local canonical coordinates $(x,\xi)$ on $T^*X$  near $(x^0,\xi^0):=\pl(c^0)$ and $(y,\eta)$ on $T^*Y$ near $(y^0,\eta^0):=\pr(c^0)$ such that, as for the model $C_0$ in
Sec. \ref{sec examps}, $(x_1,x',y_1,\eta'):=(x_1,x_2,x_3,y_1,\eta_2,\eta_3)$ form local coordinates on $C$ near $c^0$. In fact, since  $d\pl$ and $d\pr$ drop rank by 1 at $c^0$, there exist (see \cite{H-book}) symplectic decompositions $T_{(x^0,\xi^0)}(T^*X)=V_1\oplus V',\, T_{(y^0,\eta^0)}(T^*Y)=W_1\oplus W'$, with $\dim(V_1)=\dim(W_1)=2,\, \dim(V')=\dim(W')=4$, and
\[
T_{c^0}C=Gr(\chi)\oplus\big(\Lambda_L\times(0)\big)\oplus\big((0)\times\Lambda_R\big),
\]
where $\Lambda_L\subset V_1,\, \Lambda_R\subset W_1$ are Lagrangian, i.e., one-dimensional, subspaces and $\chi\in Sp(V',W')$. Thus, in suitable linear symplectic coordinates, $\Lambda_L=\{\xi_1=0\},\, \Lambda_R=\{\eta_1=0\}$ in $V_1\simeq W_1\simeq T^*\R$,  and $\chi=I\in Sp(T^*\R^2,T^*\R^2)$, so that $(x_1,y_1,x',\eta')$ are coordinates on $T_{c^0}C$.  By a standard application of the homogeneous Darboux theorem \cite{H-book, M}, there are local canonical coordinates $(x,\xi,y,\eta)$ on $T^*X\times T^*Y$ such that $c^0=(0,e_3*; 0,e_3*)$ and $(x,y_1,\eta')$ form local coordinates on $C$.

We further prepare $C$ by noting that if $v_R,v_L$ are kernel vector fields for $\pr,\pl$,  i.e., nonzero vector fields along $\Sigma_1$ generating $\ker(d\pr),\, \ker(d\pl)$, resp., then $v_L$ is a linear combination of $\partial_{y_1},\,\partial{\eta_2},\, \partial{\eta_3}$. Thus, $d\pr(v_L)$ is a vector field along $\pr(\Sigma_{1,1})$, nontangent to $\Sigma_{1,1}$ since $\pr$ is a cusp. Using homogeneous Darboux  and the nonradiality of
$\pr(\Sigma_{1,1})$, one can assume that $v_L=\partial_{y_1}$ at $\Sigma_{1,1}$, $\pr(\Sigma_{1,1})=\{\eta_1=\eta_2=0\}$ and $d\pr(TC)=span\{T\pr(\Sigma_{1,1}),d\pr(v_L)\}=\{d\eta_1=0\}$ along $\Soo$.
Working similarly on the left, one can assume $v_R=\partial_{x_1}, \, \pl(\Sigma_{1,1})=\{\xi_1=\xi_2=0\}$ and $d\pl(TC)=\{d\xi_1=0\}$ along $\Soo$.

Furthermore,  near $c^0$ one can assume that  (i) $\Sigma_1\!=\!\{\eta_2=f(x,y_1,\eta_3)\}$, and (ii)
$\Sigma_{1,1}=\{\eta_2=f(x,y_1,\eta_3),\, y_1=g(x,\eta_3)\}$ for functions $f,g$, homogeneous of degrees 1,0, resp., and with $\partial_{x_1}g\ne 0$. To see (i), let $\tf(x,y_1,\eta')$ be a defining function for $\So$, homogeneous of degree 1; by the implicit function theorem it suffices to show that $d_{\eta_2}\tf(c^0)\ne 0$. Euler's identity implies that $\eta'\cdot d_{\eta'}\tf=0$ at $\So$; since $\eta_2=0,\, \eta_3\ne\nolinebreak0$ at $c^0$, this yields $d_{\eta_3}\tf=0$.
Since $\pl,\pr$ are cusps, the kernel vector fields $v_L,v_R\in T\So$ at $\Soo$, thus at $c^0$, and hence $d_{x_1}\tf(c^0)=d_{y_1}\tf(c^0)=\nolinebreak0$. If $d_{\eta_2}\tf=0$, the only nonzero components of $d\tf(c^0)$ would be $d_{x'}\tf$, and this contradicts the fact that $d\pl(T\So)=d\pl(T\Soo)=T\pl(\Soo)=\{d\xi_1=d\xi_2=0\}$. Hence, $d_{\eta_2}\tf(c^0) \ne 0$, and thus  $\So$ can be described as $\{\eta_2=f(x,y_1,\eta_3)\}$ for some $f$. Note that since $v_R,v_L\in T\So$ at $\Soo$, one has $\partial_{x_1}f|_{\Soo}=\partial_{y_1}f|_{\Soo}=0$.
For (ii), note that since $\pl:\Soo\lra \pl(\Soo)$ is a diffeomorphism, and $(x,\xi_3)$ are coordinates on $\pl(\Soo)$, it follows that $(\pl^*(x),\pl^*(\xi_3))=(x,d_{x_3}S(x,y_1,\eta'))$ form coordinates on $\Soo$.  Since $T\So=span\{T\Soo,v\}$ for both $v=v_R=\partial_{x_1}$ and $v=v_L=\partial_{y_1}$, $(x,y_1,d_{x_3}S(x,y_1,f,\eta_3))$ are coordinates on $\So$, and we can replace $d_{x_3}S$ by $\eta_3$ since $d^2_{x_3\eta_3}S \ne 0$. Then, $\Soo$ is a hypersurface in $\So$ transverse to both $\partial_{x_1}$ and $\partial_{y_1}$, and can thus be described as $\{\eta_2=f(x,y_1,\eta_3),\, y_1=g(x,\eta_3)\}$ with $\partial_{x_1}g\ne 0$.

%\medskip

 Since $(x,y_1,\eta')$ form microlocal coordinates on $C$,  there exists a generating function $S(x,y_1,\theta')$, so that $\phi(x,y,\theta')=S(x,y_1,\theta')-y'\cdot\theta'$ parametrizes $C$ near $c^0$.
We have
\bfo
C=\Big\{ \big(x, \  \partial_xS; \   y_1,   \partial_{\theta'}S, -\partial_{y_1}S,  \theta'\big) : \, x\in\R^3,\ y_1\in\R,\theta'\in\R^2\setminus 0 \Big\},
\efo
with  $\pi_L(x, y_1, \theta')=(x, \
\partial_x S(x, y_1, \theta'))$, $\pi_R(x, y_1, \theta')=(y, \  -\partial_{y_1}S(x, y_1, \theta'), \theta')$  and $\Soo=\{\theta_2-f=y_1-g=0\}$.
Since
$\pi_L(\Sigma_{1,1})=\{\xi_1=\xi_2=0 \}$,  one has
$\partial_{x_1}S|_{\{\theta_2-f=y_1-g=0\}}=0\hbox{ and  }\partial_{x_2}S|_{\{\theta_2-f=y_1-g=0\}}=0$.
From the second relation, it follows that $S|_{\{\theta_2-f=y_1-g=0 \}}$ is independent of $x_2$,  so that  $S$ is a function of just $x_1, x_3, \theta_3$; however, from the first relation one can express $S$  as
\[
 S(x, y_1, \theta)= S_0(x_3, \theta_3)+(y_1-g)^2S_1(x, y_1, \theta') +
(\theta_2-f)S_2(x, y_1, \theta')
\]
Similarly,  $\pi_R(\Sigma_{1,1})=\{\eta_1=\eta_2=0\}$  implies  that $\partial_{y_1}S|_{\{\theta_2-f=0=y_1-g\}}=0$,
but   doesn't provide any more information about  $S$, cf. Remark 1.
On $C$, the ideal of $\Sigma_{1,1} $ equals $\left( \theta_2-f(x,y_1,\theta_3), \ y_1-g(x, \theta_3)  \right) = \left( S_{y_1}, \theta_2  \right)$, so $f= \alpha S_{y_1} + \beta \theta_2 $. Hence $f|_{\Sigma_{1,1}}=0$; $\partial_{\theta_3}f|_{\Sigma_{1,1}}=0$ (since $f$ is homogeneous of degree $1$); $f=c(y_1-g)^2$ for some $c\in C^\infty$; and  $\partial_{\theta_3}g|_{\Sigma_{1,1}}=0$, since $g$ is homogeneous of degree $0$.
Thus, the canonical relation  must have the form,
\bfo
C&=& \Big\{\big(  x_1, x_2, x_3,  -2(y_1-g)\partial_{x_1}g S_1 + (y_1-g)^2 \partial_{x_1}S_1 + (\theta_2-f) \partial_{x_1}S_2-\partial_{x_1}fS_2, \\
& &\qquad  (y_1-g)^2 \partial_{x_2}S_1 - 2(y_1-g)\partial_{x_2}gS_1 + (\theta_2-f) \partial_{x_2}S_2-\partial_{x_2}f S_2, \\
& &\qquad (y_1-g)^2 \partial_{x_3}S_1 - 2(y_1-g) \partial_{x_3}gS_1+ (\theta_2-f) \partial_{x_3} S_2 -\partial_{x_3}f S_2+ \partial_{x_3}S_0; \\
& &\qquad y_1, y_2, y_3, \ 2(y_1-g)S_1 + (y_1-g)^2 \partial_{y_1}S_1 + (\theta_2-f) \partial_{y_1}S_2-\partial_{y_1}fS_2, \\
& &\qquad -\theta_2, -\theta_3\big):
 y_2= (y_1-g)^2\partial_{\theta_2}S_1 + (\theta_2-f)\partial_{\theta_2}S_2 + S_2, \\
& &\qquad  y_3=(y_1-g)^2\partial_{\theta_3}S_1 -2(y_1-g)\partial_{\theta_3}gS_1+(\theta_2-f)\partial_{\theta_3}S_2 -\partial_{\theta_3}fS_2+  \partial_{\theta_3}S_3\Big\}.
\efo
From this, one sees  that
$$d\pi_L\equiv \left[ \begin{array} {cccccc} 1 & 0 & 0 & 0 & 0 & 0  \\
0 &  1  & 0 & 0 & 0 & 0  \\
0 & 0 & 1 & 0 & 0 & 0  \\
\cdot & \cdot & \cdot &  -2 \partial_{x_1}gS_1-\partial_{x_1y_1}fS_2 & \partial_{x_1}S_2 & 0 \\
\cdot & \cdot & \cdot &  -2 \partial_{x_2}gS_1-\partial_{x_2y_1}fS_2 & \partial_{x_2}S_2 & 0\\
\cdot & \cdot & \cdot & \cdot & \partial_{x_3}S_2 & \partial^2_{x_3 \theta_3}S_0\\
\end{array} \right] \!\!\!\!\!\mod ( y_1-g , \theta_2-f )$$
with  $\det\, d \pi_L=[(-2 \partial_{x_1}gS_1-\partial_{x_1y_1}fS_2)\partial_{x_2}S_2- ( -2 \partial_{x_2}gS_1-\partial_{x_2y_1}fS_2)\partial_{x_1}S_2] ( \partial^2_{x_3 \theta_3}S_0) \linebreak+ (y_1-g)E_1 + (\theta_2-f) E_2 $ for some $E_1,E_2$. Since $\pi_L$ has  a
cusp singularity at $\Sigma_{1,1}$, the determinant must vanish simply. We have  that $\partial^2_{x_3
\theta_3}S_0|_{\Sigma_{1,1}} \neq 0$ by homogeneity and $\partial_{x_2 \theta_2}S=\partial_{x_2}S_2|_{\Sigma_{1,1}} \neq 0$ since $d\pi_L(TC)=\{ d\xi_1=0 \}$. Since Ker d$\pi_L=\partial_{y_1}$,  the $y_1$ column in $d\pl$ must be $0$ at $\Sigma_{1,1}$. We have that $(-2 \partial_{x_1}gS_1-\partial_{x_1y_1}fS_2)|_{\Sigma_{1,1}} = \partial_{x_1}g (S_1+cS_2)|_{\Sigma_{1,1}}=0 $ so $(S_1+cS_2)|_{\Sigma_{1,1}}=0$, and
$(-2 \partial_{x_2}gS_1-\partial_{x_2y_1}fS_2)|_{\Sigma_{1,1}} = \partial_{x_2}g (S_1+cS_2)|_{\Sigma_{1,1}}=0 $, and since $\partial_{x_2}S_2, \partial_{x_1}g \ne 0$ we can conclude that det d$\pi_L|_{\Sigma_{1,1}}= (S_1+cS_2)|_{\Sigma_{1,1}}$.
Since $S_0$ is nondegenerate, by another canonical transformation on the left we may assume that $S_0=x_3\theta_3$.
From  $\ker d \pi_L =\R\cdot \frac{\partial}{\partial y_1}$, we have  $\partial_{y_1} (S_1+cS_2)|_{\Sigma_{1,1}}=0$. Since $d\pi_R(TC)=\{ d\eta_1=0 \}$ we obtain that $\partial_{y_1}S_2|_{\Sigma_{1,1}}=\partial_{\theta_2 y_1}S|_{\Sigma_{1,1}}=0$.
One also has that $\partial_{y_1}S_1|_{\Sigma_{1,1}}=0$, since a calculation shows that this is $\partial^2_{y_1}\eta_1(x,y_1,\theta')|_{\Soo}$, which vanishes due to the cusp structure of  $d\pr$ and the fact that
$d\pr(TC)=\{d\eta_1=0\}$.  Then one can write $S_1=(y_1-g)^2 S_5 + (\theta_2-f) S_6$ with $S_5 \neq 0$,
leading to  a weak normal form for a phase function parametrizing the flat two-sided cusp $C$, namely
$\phi(x, y, \theta_2,\tht) = \left(S_2 + (y_1-g)^2S_6 \right)(\theta_2-f)+\linebreak(x_3-y_3)\theta_3 + (y_1-g)^4 S_5 -y_2\theta_2$,
with $\partial_{x_2} S_2(x,y_1,\theta')|_{\Sigma_{1,1}} \neq 0, \quad S_5, S_6 \neq 0$.

\par In studying the composition in the next section, $f,\, g$  introduce only algebraic complications and do not affect the final result.  Thus, for simplicity we take $f=0$ and $g=x_1$ and will use the phase function
 \beq\label{eqn wnf}
\phi(x, y, \theta_2,\tht)& =& \left(S_2-y_2 + (x_1-y_1)^2S_4\right)\theta_2+(x_3-y_3)\theta_3 + (x_1-y_1)^4 S_3 ,\nonumber \\
& &\qquad\qquad\qquad \partial_{x_2} S_2(x,y_1,\theta') \neq 0, \quad S_3, S_4 \neq 0 \hbox{ on } \Sigma_{1,1},
\eeq
with $S_5, S_6$ relabeled as $S_3, S_4$, resp.

\section{Composition}\label{sec comp}

\par

Now consider the composition $B^*A$ for  $A\in I^m(X,Y;C),\, B\in I^{m'}(X,Y;C)$.
Conjugating by unitary FIOs associated with the canonical transformations used in the last section, we may assume that $C$ is parametrized by the phase function $\phi$ from  (\ref{eqn wnf}). We have, for $a\in S^{m+\frac12},\, b\in S^{m'+\frac12}$,
\bfo
B^*A f(x)=\int e^{i(\phi(z, y, \eta')-\phi(z, x, \theta'))} a (z, y, \eta') \bar{b}(z, x, \theta') f(y) dy d z d \theta' d \eta',
\efo
where the phase $\phi(z, y, \eta')-\phi(z, x, \theta')$, which we denote $\hat\phi$, is
\bfo
 (z_3-y_3) \eta_3 &+& (z_1-y_1)^4 S_3(z, y_1, \eta') +\left(S_2(z, y_1,
\eta')-y_2+(z_1-y_1)^2 S_4 (z, y_1, \eta')\right)\eta_2 \\
-  (z_3-x_3)\theta_3&-&(z_1-x_1)^4 S_3(z, x_1, \theta')  -\left(S_2(z, x_1,
\theta')-x_2+(z_1-x_1)^2 S_4(z, x_1,\theta')\right)\theta_2.
\efo
We perform  stationary phase in $z_2, \eta_2$ and $z_3, \eta_3$:
\beq\label{eqn ztwo}
\partial_{z_2}\hat{\phi} &\!\!\!=\!\!\!&  \eta_2\partial_{z_2}S_2- \theta_2\partial_{z_2}S_2  + \partial_{z_2}S_3 (z_1-y_1)^4 -\partial_{z_2}S_3 (z_1-x_1)^4 \\
& & + (z_1-y_1)^2\eta_2 \partial_{z_2}S_4  -(z_1-x_1)^2\theta_2\partial_{z_2}S_4, \nonumber \\
\label{eqn zthree}
\partial_{z_3} \hat{\phi}&\!\!\!=\!\!\!& \eta_3- \theta_3 + (z_1-y_1)^4 \partial_{z_3} S_3 -(z_1-x_1)^4 \partial_{z_3} S_3 +
\eta_2\partial_{z_3}S_2  -\theta_2 \partial_{z_3} S_2\nonumber  \\
& &  +(z_1-y_1)^2  \eta_2 \partial_{z_3} S_4-(z_1-x_1)^2
\theta_2\partial_{z_3}S_4, \\
\label{eqn etatwo}
\partial_{\eta_2} \hat{\phi}&\!\!\!=\!\!\!& (z_1-y_1)^4 \partial_{\eta_2}S_3 + S_2-y_2 +(z_1-y_1)^2[S_4+  \eta_2\partial_{\eta_2} S_4  ]+ \eta_2\partial_{\eta_2}S_2 , \\
\label{eqn etathree}
\partial_{\eta_3} \hat{\phi}&\!\!\!=\!\!\!& z_3-y_3 + (z_1-y_1)^4 \partial_{\eta_3}S_3 +  \eta_2\partial_{\eta_3}S_2 + (z_1-y_1)^2\eta_2
\partial_{\eta_3}S_4.
\eeq
 We solve for $ z_3$ using (\ref{eqn etathree}) and find
\bfo
z_3&=&y_3 -(z_1-y_1)^4 \partial_{\eta_3}S_3 -\eta_2\partial_{\eta_3}S_2  - (z_1-y_1)^2 \eta_2\partial_{\eta_3} S_4 := y_3+S_5.
\efo
From  equation (\ref{eqn etatwo}) we can  solve  for $z_2$ since $\partial_{z_2}S_2 \neq 0$, but only implicitly.
Now rewrite (\ref{eqn zthree}) as
\bfo
\eta_3&=&\theta_3 + (z_1-x_1)^4 \partial_{z_3}S_3 -(z_1-y_1)^4\partial_{z_3}S_3\\
& &\quad \, + (z_1-x_1)^2 \theta_2\partial_{z_3}S_4
-(z_1-y_1)^2 \eta_2 \partial_{z_3}S_4+ \theta_2\partial_{z_3}S_2 -\eta_2\partial_{z_3}S_2,
\efo
and solve for $\eta_3$, using the fact that $\theta_2\partial_{z_3}S_2 -\eta_2\partial_{z_3}S_2 $ vanishes at $\{x_1=y_1,\eta'=\theta'\}$, so that
$$\theta_2\partial_{z_3}S_2 - \eta_2\partial_{z_3}S_2\!=\! (x_1-y_1) \theta_2  \partial^2_{z_3 x_1}S_2
+(\theta_2-\eta_2)(\partial_{z_3}S_2 +\theta_2 \partial^2_{z_3 \theta_2}S_2 ) + (\theta_3-\eta_3)\theta_2 \partial^2_{z_3 \theta_3}S_2 .$$
The other differences in the equations (\ref{eqn ztwo}),(\ref{eqn etathree})
 can be handled  similarly, so
$$\eta_3=\theta_3 + (x_1-y_1)M_1 (z_1, z_2, x_1, y_1, \theta', \eta_2) + (\theta_2-\eta_2)N_1(z_1, z_2, x_1, y_1, \theta', \eta_2) + $$ $$(\theta_3-\eta_3) P(z_1, z_2, x_1, y_1, \theta', \eta_2)$$
for some functions $M_1,N_1$ and $P$.
Solving for $\eta_3-\theta_3$, we have $$\eta_3-\theta_3=[(x_1-y_1)M_1 +(\theta_2-\eta_2)N_1 ](I+P)^{-1}= (x_1-y_1)M_2 +
(\theta_2-\eta_2)N_2.$$
Also,
$\eta_2-\theta_2= (x_1-y_1)M_3 + (\eta_3-\theta_3)N_3$, so that $\eta_3-\theta_3=(x_1-y_1)M_4$
 and $\eta_2-\theta_2=(x_1-y_1)N_4$.
Hence, the phase becomes

\bfo
\hat{\phi}&=& (x_3-y_3)\theta_3 + (x_2-y_2)\theta_2 + \left((z_1-y_1)^2S_4(z_1, y,
\theta')-(z_1-x_1)^2S_4(z_1, x_1, y_2, y_3, \theta')\right)\theta_2 \\
& &+ (z_1-y_1)^4 S_3-(z_1-x_1)^4S_3 + \left(S_2(z_1, y, \theta')-S_2(z_1, x_1, y_2, y_3, \theta')\right)\theta_2 + (x_1-y_1)S_6.
\efo
We have
\bfo
\big((z_1-y_1)^2S_4(z_1, y, \theta')&-&(z_1-x_1)^2S_4(z_1, x_1, y_2, y_3, \theta')\big)\theta_2 \\
 &=& [(z_1-y_1)^2S_4(z_1, y_1, \cdot)-(z_1-x_1)^2S_4(z_1, y_1, \cdot)]\theta_2\\
 & &+ [(z_1-x_1)^2 S_4(z_1, y_1, \cdot)-(z_1-x_1)^2S_4(z_1, x_1, \cdot)]\theta_2\\
 &= & [(x_1-y_1)(2z_1-x_1-y_1 )S_4(z_1,y_1,\cdot) + (z_1-x_1)^2 (x_1-y_1) \partial_{y_1}S_4]\theta_2,
 \efo

 \bfo
 (z_1-y_1)^4 S_3-(z_1-x_1)^4S_3&=&(z_1-y_1)^4 S_3(z_1, y_1, \cdot) -(z_1-x_1)^4 S_3(z_1, y_1, \cdot) \\
 && + (z_1-x_1)^4S_3(z_1, y_1, \cdot)- (z_1-x_1)^4S_3(z_1, x_1,\cdot) \\
 &=& (x_1-y_1)(2z_1-x_1-y_1)((z_1-y_1)^2+(z_1-x_1)^2)S_3(z_1,y_1,\cdot)\\
 & & +(z_1-x_1)^4 (x_1-y_1)\partial_{y_1}S_3,
\efo
 and can write
$(z_1-x_1)^2 (x_1-y_1) \partial_{y_1}S_4\theta_2+(z_1-x_1)^4 (x_1-y_1)\partial_{y_1}S_3
 +\linebreak\left(S_2(z_1, y, \theta')-S_2(z_1, x_1, y_2, y_3, \theta')\right)\theta_2 + (x_1-y_1)S_6:=(x_1-y_1)S_7$.
From this, it follows that $\hat{\phi}=(x_2-y_2)\theta_2 + (x_3-y_3)\theta_3 + (x_1-y_1)\big[(2z_1-x_1-y_1)((z_1-y_1)^2+(z_1-x_1)^2)S_3   + (2z_1-x_1-y_1)S_4\theta_2 + S_7(z_1, y, x_1, \theta)\big]$.
 We make a change of variables similar to that in the model case in Sec. \ref{subsec comp}:
 $$\theta_1= (2z_1-x_1-y_1)((z_1-y_1)^2+(z_1-x_1)^2)S_3 +(2z_1-x_1-y_1)S_4\theta_2 +S_7, $$
 \bfo
 \frac{\partial \theta_1}{\partial z_1}&=&2S_4 \theta_2+(2z_1-x_1-y_1)\theta_2\partial_{z_1}S_4 +
 2S_3((z_1-x_1)^2+(z_1-y_1)^2) \\
 & &\quad+2(2z_1-x_1-y_1)^2S_3+ (2z_1-x_1-y_1)((z_1-y_1)^2+(z_1-x_1)^2)\partial_{z_1}S_3 + \partial_{z_1}S_7
 \efo
and
\bfo
\frac{\partial^2 \theta_1}{\partial z_1^2}&=& 4 \partial_{z_1}S_4 \theta_2 + (2z_1-x_1-y_1)\theta_2\partial^2_{z_1^2}S_4 +  4
((z_1-x_1)^2 +(z_1-y_1)^2)\partial_{z_1}S_3 \\
&& +12 (2z_1-x_1-y_1)S_3
+ 4(2z_1-x_1-y_1)^2 \partial_{z_1}S_3 +
(2z_1-x_1-y_1)((z_1-y_1)^2\\
&&+(z_1-x_1)^2)\partial^2_{z_1}S_3 + \partial_{z_1}S_7.
\efo
Near points on $\Sigma_{1,1}$,  $\frac{\partial^2 \theta_1}{\partial z_1^2}=0$ implies $z_1=\frac{x_1+y_1}{2} + R(x_1,y,
\theta_2, \theta_3)$, for some $R$ which vanishes at the base point $(x_1,y, \theta_2, \theta_3) = (0,0, 0,0,0,1)$. Then
\bfo
\frac{\partial \theta_1}{\partial z_1}&=&2S_4 \theta_2+ R \partial_{z_1}S_4 \theta_2 + 2S_3(\frac{(x_1-y_1)^2}{2}+2R^2) \\
& & +2R^2S_3+ R(\frac{(x_1-y_1)^2}{2}+2R^2)\partial_{z_1}S_3 + \partial_{z_1}S_7=0
\efo
implies that
$$\theta_2\!=\!-\frac{(x_1-y_1)^2}{2}\frac{2S_3+R\partial_{z_1}S_3}{2S_4+R\partial_{z_1}S_4}+P_1(x_1,y
,\theta_3)\!=\!-\frac{(x_1-y_1)^2}{2}N(x_1,y,\theta_3)+P_1(x_1,y ,\theta_3)$$
for some $N,P_1$,
and
$$\theta_1=R(\frac{(x_1+y_1)^2}{2}+2R^2)S_3+RS_4\theta_2+S_7:=  P_2(x_1,y, \theta_3).$$
Since $S_3, S_4 \neq 0$ and $R=0$ at the base point, one has $N \neq 0$.
 Following
  the same analysis as in the model case, one sees that
  $$\hat{\phi}= (x_2-y_2)\theta_2 + (x_3-y_3)\theta_3 + (x_1-y_1)\theta_1+
 (\frac{\theta_1-P_2}{\theta_3}-t_1t_2-\frac{1}{2}t_1^3)\sigma_1 +$$
 $$(\frac{\theta_2}{\theta_3}+\frac{(x_1-y_1)^2}{2}\frac{N}{\theta_3}-\frac{P_1}{\theta_3}-t_2)\sigma_2.$$
\par Next, we perform a stationary phase in $t_2, \sigma_2$:
\[
\partial_{t_2} \hat{\phi}=-t_1 \sigma_1-\sigma_2,\quad \partial_{\sigma_2}\hat{\phi}=\frac{\theta_2}{\theta_3}+\frac{N}{
2\theta_3}(x_1-y_1)^2-\frac{P_1}{\theta_3}-t_2
\]
and
$$\hat{\phi}=(x-y)\theta+ (\frac{\theta_1}{\theta_3}  -\frac{P_2}{\theta_3}-t_1[\frac{\theta_2}{\theta_3}+\frac{N}{2
\theta_3} (x_1-y_1)^2-\frac{P_1}{\theta_3}]-\frac{1}{2}t_1^3)\sigma_1   $$
\par With one more stationary phase in $\sigma_1$ and $\theta_1$, obtaining
\bfo
\partial_{\theta_1}\hat{\phi}= x_1-y_1+\frac{\sigma_1}{\theta_3},\quad
\partial_{\sigma_1} \hat{\phi}= \frac{\theta_1}{\theta_3}  -\frac{P_2}{\theta_3}-t_1[\frac{\theta_2}{\theta_3}+\frac{N}{2
\theta_3} (x_1-y_1)^2-\frac{P_1}{\theta_3}] -\frac{1}{2}t_1^3,
\efo
and then homogenize $t_1$ by setting
$\tau=t_1 \theta_3$. The resulting phase is then
\beq\label{eqn star}
\quad\hat{\phi}=\left(x_2-y_2+\frac{\tau}{\theta_3}(x_1-y_1)\right)\theta_2 &+&\Big(x_3-y_3+(x_1-y_1)
\frac{P_2-\frac{\tau}{\theta_3}P_1}{\theta_3}\\
&+& (x_1-y_1)^3\frac{\tau}{\theta_3^2}\frac{N}{2}+ (x_1-y_1)\frac{1}{2}(\frac{\tau}{\theta_3})^3\Big)\theta_3,\nonumber
\eeq
and the Schwartz kernel has the representation
\beq\label{eqn bstara}
K_{F^*F}(x,y)=\int e^{i \hat{\phi}} c(x,y,\theta',\tau)\, d\theta_2 d\theta_3 d \tau, \quad c\in S^{m+m'}.
\eeq
For each of the variables $s=x_1,y_1,y_2,y_3$, let $[[s]]$ denote
\[
[[s]]\!\!=\!
\!\partial_{s}\big({P_2-\frac{\tau}{\theta_3}P_1}\big)+(x_1-y_1)^2\frac{\tau}{\theta_3}
\partial_{s}\frac{N}{2},\,
[[\theta_3]]\!=\!\partial_{\theta_3}\big(P_2-\frac{\tau}{\theta_3}P_1+
(x_1-y_1)^2\frac{\tau}{\theta_3}\frac{N}{2}+\frac{1}{2}\frac{\tau^3}{\theta_3^2}\big).
\]
Then we have
\bfo
WF(F^*F)\!\!\!\!&\subset &\!\!\!\!\Big\{\big(x_1, x_2, x_3, \frac{\tau \theta_2}{\theta_3}+3(x_1-y_1)^2
\tau \frac{N}{2\theta_3} + P_2-\frac{\tau}{\theta_3}P_1+ \frac{1}{2}(\frac{\tau}{\theta_3})^3 \theta_3+ (x_1-y_1)[[x_1]] ,
\theta'; \\
& &y_1, y_2, y_3,\frac{\tau \theta_2}{\theta_3}+3(x_1-y_1)^2
\tau \frac{N}{2\theta_3} + P_2-\frac{\tau}{\theta_3}P_1+ \frac{1}{2}(\frac{\tau}{\theta_3})^3\theta_3+ (x_1-y_1)[[y_1]],\\
& & \theta_2 + (x_1-y_1)[[y_2]], \theta_3 + (x_1-y_1)[[y_3]]\big) :\\
& & \quad d_{\tau}\hat{\phi}=(x_1-y_1)\left(\frac{\theta_2-P_1}{\theta_3}+\frac{3}{2} (\frac{\tau}{\theta_3})^2
+(x_1-y_1)^2\frac{N}{2\theta_3 }\right)=0,\\
& &\quad  d_{\theta_2}\hat{\phi}=x_2-y_2+(x_1-y_1)\frac{\tau}{\theta_3}=0,\\
& &\quad d_{\theta_3}\hat{\phi}=x_3-y_3-(x_1-y_1)\frac{\tau \theta_2}{\theta_3^2}+(x_1-y_1)[[\theta_3]] =0\Big \}.
\efo
The phase function is degenerate because the equation $d_{\tau}\hat{\phi}=0$ has a normal crossing; from $\{x_1-y_1=0\}$,  the contribution to $WF(B^*A)$ is contained in
\bfo
\Delta&=&\Big\{ \big(x, \frac{\tau \theta_2}{\theta_3}+ P_2-\frac{\tau}{\theta_3}P_1+
\frac{1}{2}(\frac{\tau}{\theta_3})^3\theta_3, \theta_2, \theta_3; \\
& &\quad x,\frac{\tau \theta_2}{\theta_3}+ P_2-\frac{\tau}{\theta_3}P_1+ \frac{1}{2}(\frac{\tau}{\theta_3})^3\theta_3,
\theta_2, \theta_3\big): x\in\R^3,(\tau,\theta')\in\R^3\setminus 0 \Big \},
\efo
which is a cusp parametrization of  $\Delta$ (as in the model case). On the other hand, if $x_1-y_1 \neq 0$ then

\bfo
\theta_2&=&P_1-\frac{3}{2}(\frac{\tau}{\theta_3})^2 \theta_3-(x_1-y_1)^2\frac{N}{2},\quad
y_2=x_2+(x_1-y_1)\frac{\tau}{\theta_3}, \\
y_3&=&x_3-(x_1-y_1)\frac{\tau \theta_2}{\theta_3^2}+(x_1-y_1)[[\theta_3]], \\
\xi_1&=&P_2- (\frac{\tau}{\theta_3})^3\theta_3+\tau N\frac{(x_1-y_1)^2}{\theta_3} +(x_1-y_1)[[x_1]], \\
\eta_1&=&P_2- (\frac{\tau}{\theta_3})^3\theta_3+\tau N\frac{(x_1-y_1)^2}{\theta_3} +(x_1-y_1)[[y_1]],  \\
\eta_2&=&\theta_2+(x_1-y_1)[[y_2]], \hbox{ and }
\eta_3=\theta_3+(x_1-y_1)[[y_3]].
\efo
Let $\tilde{C}:=\tilde{C}_{\hat{\phi}}$ be the image of $\psi(x_1, x_2, x_3, y_1, \theta_3, \tau)=(x,\xi;y,\eta) $ with $\xi_2=\theta_2,
\xi_3=\theta_3$, and $\xi_1, \eta_1, \eta_2, \eta_3, y_2, y_3$ given above. As in the model case, $\psi$ satisfies
  Def. \nolinebreak\ref{def ou}: it drops rank by one simply  at  $\Sigma_1=\{x_1-y_1=0=\tau\}$, and  $\ker
d\pi=\R\cdot\frac{\partial}{\partial \tau} \nsubseteq T\Sigma$, so that $\tilde{C}$ is an open umbrella.
%\medskip

\subsection{Generalized Fourier integral operators for $\Delta\cup \tilde{C}$}\label{sec gfios}

Using the phase $\hat\phi$, one can now define a  class of generalized FIOs,  $G$, with the $WF(G)$ contained in the union of the diagonal and the open umbrella. For a general $\tilde{a} \in S^{\mu}$,
define  $G$ by the right hand side of (\ref{eqn bstara}), but with the amplitude $\tilde{a}$;
then one has $WF(G) \subset \Delta \cup \tilde{C}$.
Furthermore, one can compute the orders of $G$ on $\Delta \setminus \tilde{C}$ and $\tilde{C} \setminus  \Delta$,  as in the model case at the end of Sec. \ref{subsec comp}.
Due to the normal crossing of $d_\tau\hat\phi$,  the critical set $Crit_{\hat{\phi}}= \{ (x, y, \tau, \theta'): d_{\tau} \hat{\phi}=0, d_{\theta'} \hat{\phi}=0 \}$ decomposes as $Crit_{\Delta} \cup Crit_{\tilde{C}} $, and $\sigma=\rho (E_{\hat{\phi}})^{\frac{1}{2}}$, where $E_{\hat{\phi}}=|\frac{D(\lambda_i, \frac{\partial
 \hat{\phi}}{\partial \theta})}{D(x, \theta)}|^{-1}$ and $\lambda_i$ are  local coordinates  on $Crit_{\hat{\phi}}$.

  \par On $Crit_{\Delta}$, local coordinates are $(x,\tau, \theta')$ and
  \bfo
  E_{\hat{\phi}}=\Big|\frac{D(x,\tau,
  \theta', \frac{\partial \hat{\phi}}{\partial \tau}, \frac{\partial \hat{\phi}}{\partial \theta_2},
  \frac{\partial \hat{\phi}}{\partial \theta_3})}{D(x, \tau, \theta', y_1, y_2, y_3)}\Big|^{-1}=\Big(
  \frac{\theta_2}{\theta_3} + \frac{3}{2}\frac{\tau^2}{\theta_3^2}-\frac{P_1}{\theta_3}\Big)^{-1};
  \efo
  hence  $G \in I^{\mu} (\Delta \setminus \tilde{C})$ and the principal symbol satisfies $\sigma \sim \delta^{-\frac{1}{2}}$, where $\delta$ is the distance to $\Delta \cap \tilde{C}$.
On the other hand, $(x,y_1, \tau, \theta_3)$ are local coordinates  on $Crit_{C_{\hat{\phi}}}$,
 \bfo
 E_{\hat{\phi}}&=&\Big|\frac{D(x, y_1, \tau,
  \theta_3 , \frac{\partial \hat{\phi}}{\partial \tau}, \frac{\partial \hat{\phi}}{\partial \theta_2}, \frac{\partial
  \hat{\phi}}{\partial \theta_3})}{D(x, y_1, \tau, \theta_2, \theta_3, y_2, y_3, \theta_2)}\Big|^{-1} \\
  &=&\Big(\frac{x_1-y_1}{\theta_3}+\frac{(x_1-y_1)^2}{\theta_3}[\partial_{y_3 \theta_3}P_2-\frac{\tau}{\theta_3} \partial_{y_3 \theta_3}P_1 + (x_1-y_1)^2 \frac{\tau}{2 \theta_3} \partial_{y_3 \theta_3}N] \Big)^{-1} \\
  &  \sim& \big(
  \frac{x_1-y_1}{\theta_3}\big)^{-1},
  \efo
$G \in I^{\mu} ( \tilde{C} \setminus \Delta)$
  and again $\sigma \sim\delta^{-\frac{1}{2}}$.
 The equality of the orders on  $\Delta$ and $\tilde{C}$ away from $\Delta \cap \tilde{C}$ is consistent with the composition result in \cite{N,F} for two-sided folds. This gives a precise description of the nonremovable artifact in the linearized seismic inversion problem for the single source geometry in the presence of cusp caustics.

  \section*{Acknowledgements}
  We would like to thank Cliff Nolan for patiently explaining some of  the calculations in \cite{N}. The second author
was  supported by NSF grants DMS-0551894 and -0853892.

\vskip.2in

\noindent{\sc School of Mathematical Sciences}

\noindent{\sc Rochester  Institute of Technology}

\noindent{\sc  Rochester, NY 14623}

\noindent{\tt{rxfsma@rit.edu}}

\vskip.2in

\noindent{\sc Department of Mathematics}

\noindent{\sc University of Rochester}

\noindent{\sc Rochester, NY 14627}

\noindent{\tt{allan@math.rochester.edu}}

\end{document}